\documentclass[11pt,a4paper]{amsart}
\usepackage{graphicx,multirow,array,amsmath,amssymb,color,adjustbox,kotex}

\newcommand{\veq}{\mathrel{\rotatebox{90}{$=$}}}
\def\stackbelow#1#2{\underset{\displaystyle\overset{\displaystyle\veq}{#2}}{#1}}
\newcommand{\adj}{\! \thicksim \!}
\newcommand{\nadj}{\! \nsim \!}

\newtheorem{theorem}{Theorem}

\newtheorem{proposition}[theorem]{Proposition}
\newtheorem{corollary}[theorem]{Corollary}

\title{Bipartite intrinsically knotted graphs with 23 edges}

\author[H. Kim]{Hyoungjun Kim}
\address{College of General Education, Kookmin University, Seoul 02707, Korea}
\email{kimhjun@kookmin.ac.kr}
\author[T. Mattman]{Thomas Mattman}
\address{Department of Mathematics and Statistics, California State University, Chico, Chico CA 95929-0525, USA}
\email{TMattman@CSUChico.edu}
\author[S. Oh]{Seungsang Oh}
\address{Department of Mathematics, Korea University, Seoul 02841, Korea}
\email{seungsang@korea.ac.kr}

\thanks{The first author was supported by the National Research Foundation of Korea (NRF) grant funded by the Korea government Ministry of Science and ICT(NRF-2018R1C1B6006692).}
\thanks{The third author was supported by the National Research Foundation of Korea(NRF) grant funded by the Korea government(MSIP) (No. NRF-2017R1A2B2007216).}

\begin{document}

\maketitle

\begin{abstract}
A graph is intrinsically knotted if every embedding contains a nontrivially knotted cycle.
It is known that intrinsically knotted graphs have at least 21 edges 
and that there are exactly 14 intrinsically knotted graphs with 21 edges,
in which the Heawood graph is the only bipartite graph.
The authors showed that 
there are exactly two graphs with at most 22 edges that are minor minimal bipartite intrinsically knotted: 
the Heawood graph and Cousin 110 of the $E_9+e$ family.
In this paper we show that 
there are exactly six bipartite intrinsically knotted graphs with 23 edges so that every vertex has degree 3 or more. 
Four among them contain the Heawood graph and 
the other two contain Cousin 110 of the $E_9+e$ family.
Consequently, there is no minor minimal intrinsically knotted graph with 23 edges that is bipartite.
\end{abstract}

\section{Introduction}\label{sec:int}

A graph is {\it intrinsically knotted} if every embedding of the graph in $\mathbb{R}^3$ contains a non-trivially knotted cycle.
We say graph $H$ is a {\it minor} of graph $G$ if $H$ can be obtained from a subgraph of $G$ by contracting edges.
A graph $G$ is {\it minor minimal intrinsically knotted} if $G$ is intrinsically knotted but no proper minor is.
Robertson and Seymour's~\cite{RS} Graph Minor Theorem implies that there are only finitely many minor minimal intrinsically knotted graphs.
While finding the complete list of minor minimal intrinsically knotted graphs remains an open problem, there has been recent progress
in understanding the condition for small size.

\begin{figure}[h!]
\includegraphics[scale=1]{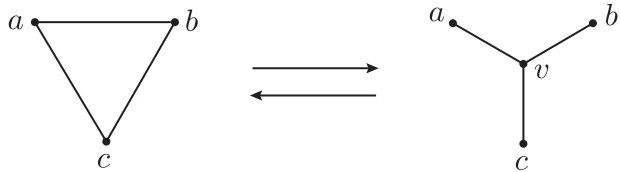}
\caption{$\nabla Y$ and $Y \nabla$ moves}
\label{fig:delta}
\end{figure}

The known examples mainly belong to $\nabla Y$ families.
A {\it $\nabla Y$ move} is an exchange operation on a graph that removes all edges of a 3-cycle $abc$ and then adds a new vertex $v$ that is connected to each vertex of the 3-cycle, as shown in Figure~\ref{fig:delta}. 
The reverse operation is a  {\it $Y \nabla$ move}.
We say two graphs $G$ and $G'$ are {\it cousins} if $G'$ is obtained from $G$ by a finite sequence of $\nabla Y$ and $Y \nabla$ moves.
The set of all cousins of $G$ is called the {\it $G$ family}.

\begin{figure}[h!]
\includegraphics[scale=1]{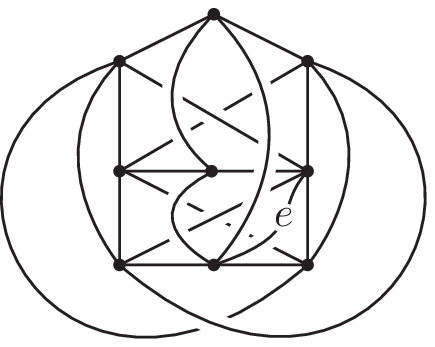}
\caption{$E_9+e$}
\label{fig:e9}
\end{figure}

Johnson, Kidwell and Michael~\cite{JKM}, and, independently, Mattman~\cite{M}, showed that intrinsically knotted graphs have at least 21 
edges. 
Lee, Kim, Lee and Oh~\cite{LKLO}, and, independently, Barsotti and Mattman~\cite{BM} showed that the complete set of minor minimal
intrinsically knotted graphs with 21 edges consists of fourteen graphs: $K_7$ and the 13 graphs obtained from $K_7$ by $\nabla Y$ moves.
There are 92 known examples of size 22: 58 in the $K_{3,3,1,1}$ family, 33 in the $E_9+e$ family (Figure~\ref{fig:e9}) and a $4$--regular
example due to Schwartz (see \cite{FMMNN}). 
We are in the process of determining whether or not this is a complete list \cite{KLLMO, KMO2, KMO3}.

In the current article, we continue a study of intrinsic knotting of bipartite graphs of small size. 
A {\em bipartite\/} graph is a graph whose vertices can be divided into two disjoint sets $A$ and $B$ such that every edge connects a vertex in $A$ to one in $B$. 
In an earlier paper~\cite{KMO}, we proved that the Heawood graph (of size 21) is the only bipartite graph among the minor minimal intrinsically knotted graphs of size 22 or less.
We also showed that Cousin 110 of the $E_9+e$ family (of size 22) is the only other graph of 22 or fewer edges that is
bipartite and intrinsically knotted and has no proper minor with both properties.
We can think of Cousin 110 as being constructed from $K_{5,5}$ through deletion of the edges in a $3$-path.
In the current paper, we extend the classification to graphs of size 23.

\begin{theorem} \label{thm:main}
There are exactly six bipartite intrinsically knotted graphs with 23 edges so that every vertex has degree 3 or more. Two of these are obtained from Cousin 110 of the $E_9+e$ family by adding an edge, the other four from the Heawood graph by adding 2 edges.
\end{theorem}

The two graphs obtained from Cousin 110 of the $E_9+e$ family are described in Subsection~\ref{subsec:110}.
Three of the graphs obtained from the Heawood graph are found in Subsection~\ref{subsec:89}
and the last in Subsection~\ref{subsec:106025}.
Since a minor minimal intrinsically knotted graph must have 
minimum degree at least three we have the following.

\begin{corollary}
There is no minor minimal intrinsically knotted graph with 23 edges that is bipartite.
\end{corollary}

The remainder of this paper is a proof of Theorem~\ref{thm:main}. 
In the next section we introduce some terminology.  Section~\ref{sec:restoring} reviews the restoring method
and introduces the twin restoring method.
Section~\ref{sec:g6} treats the case of a vertex of degree 6 or more, 
Section~\ref{sec:ab5}, the case where both $A$ and $B$ have degree 5 vertices,
and Section~\ref{sec:a5}, the case where only $A$ has degree 5 vertices.
Finally, we conclude the argument with Section~\ref{sec:g4}, which deals with the remaining cases.

\section{Terminology and strategy} \label{sec:term}

We use notation and terminology similar to that of our previous paper~\cite{KMO}.
Let $G = (A,B,E)$ denote a bipartite graph with 23 edges 
whose partition has the parts $A$ and $B$ with $E$ denoting the edges of the graph.
For distinct vertices $a$ and $b$, 
let $G \setminus \{ a,b \}$ denote the graph obtained from $G$ by deleting the two vertices $a$ and $b$.
Deleting a vertex means removing the vertex, interiors of all edges adjacent to the vertex 
and remaining isolated vertices.
Let $G_{a,b}$ denote the graph obtained from $G \setminus \{ a,b \}$ by deleting all degree 1 vertices, 
and $\widehat{G}_{a,b}=(\widehat{V}_{a,b}, \widehat{E}_{a,b})$ denote the graph obtained from $G_{a,b}$ 
by contracting edges adjacent to degree 2 vertices, 
one by one repeatedly, until no degree 2 vertex remains.
The degree of $a$, denoted by $\deg(a)$, is the number of edges adjacent to $a$.
We say that $a$ is adjacent to $b$, denoted by $a \adj b$, if there is an edge connecting them.
If they are not adjacent, we write $a \nadj b$.
If $a$ is adjacent to vertices $b, \dots, b'$, then we write $a \adj \{b, \dots, b'\}$.
If each of $a, \dots, a'$ is adjacent to all of $b, \dots, b'$, 
then we similarly write $\{a, \dots, a'\} \adj \{b, \dots, b'\}$.
Note that $\sum_{a \in A} \deg(a) = \sum_{b \in B} \deg(b) = 23$ by the definition of bipartition.
We need some notation to count the number of edges of $\widehat{G}_{a,b}$.
\begin{itemize}
\item $V_n(a)=\{c \in V\ |\ a \adj c,\ \deg(c)=n \}$
\item $V_n(a,b)=V_n(a) \cap V_n(b)$
\item $E(a)=\{ e \in E\ |\ e \,\, {\rm is \,\, adjacent \,\, to}\,\, a \}$
\item $G \setminus \{a,b\}$ has $NE(a,b)=|E(a)\cup E(b)|$ fewer edges than $G$
\end{itemize}

\begin{figure}[h]
\includegraphics{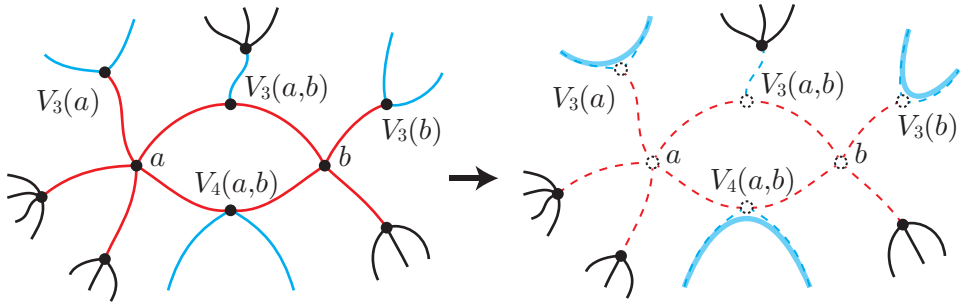}
\caption{Deriving $\widehat{G}_{a,b}$}
\label{fig21}
\end{figure}

Furthermore $G \setminus \{a,b\}$ has degree 1 or 2 vertices from $V_3(a) \cup V_3(b)$ 
and degree 2 vertices from $V_4(a,b)$ as shown in Figure~\ref{fig21}.
To derive $\widehat{G}_{a,b}$, we delete and contract the edges related to these vertices.
The total number of these edges is the sum of the following two values: 
\begin{itemize}
\item $NV_3(a,b)=|V_3(a) \cup V_3(b)| = |V_3(a)|+|V_3(b)|-|V_3(a,b)|$ 
\item $NV_4(a,b)=|V_4(a,b)|$ 
\end{itemize}

\begin{figure}[h]
\includegraphics{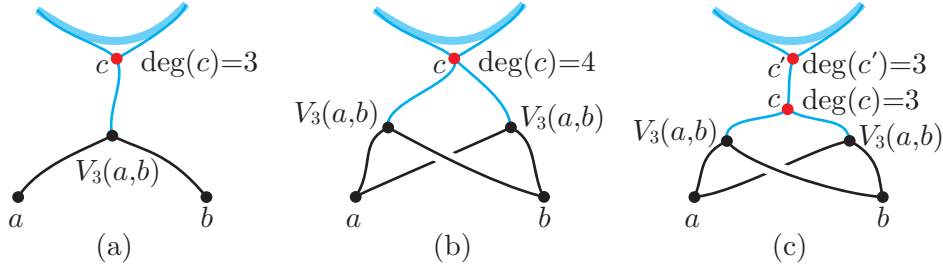}
\caption{Three types of vertex in $V_Y(a,b)$}
\label{fig22}
\end{figure}

To count $|\widehat{E}_{a,b}|$ more precisely, we need to consider the following set.
\begin{itemize}
\item $V_Y(a,b)$ is the set of removed vertices to derive $\widehat{G}_{a,b}$, 
that are adjacent to neither $a$ nor $b$; let $NV_Y(a,b)=|V_Y(a,b)|$.
\end{itemize}
This vertex set has three types as illustrated in Figure~\ref{fig22}.
In the figure, a vertex $c$ of $V_Y(a,b)$ has degree 1 or 2, 
and so must be removed in $\widehat{G}_{a,b}$.
Especially, in the rightmost figure, the two vertices $c$ and $c'$ of $V_Y(a,b)$ are removed.

Combining these ideas, we have the following equation for the number of edges of $\widehat{G}_{a,b}$,
which is called the {\em count equation\/}:
$$ |\widehat{E}_{a,b}| = 23 - NE(a,b) - NV_3(a,b) - NV_4(a,b) - NV_Y(a,b). $$

\begin{figure}[h]
\includegraphics{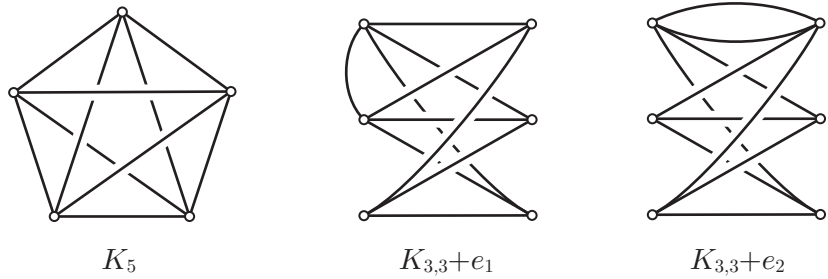}
\caption{Three non-planar graphs with 10 edges}
\label{fig:nonpla}
\end{figure}

A graph is called 2-{\em apex\/} if it can be made planar by deleting two or fewer vertices.
It is known that if $G$ is 2-apex, then it is not intrinsically knotted~\cite{BBFFHL, OT}.
So we check whether or not $\widehat{G}_{a,b}$ is planar.
The unique non-planar graph with nine edges is $K_{3,3}$.
For non-planar graphs with 10 edges, 
we consider which graphs could be isomorphic to $\widehat{G}_{a,b}$.
Note that $\widehat{G}_{a,b}$ consists of vertices with degree larger than 2.
Furthermore, $\widehat{G}_{a,b}$ may have multiple edges.
There are exactly three non-planar graphs on 10 edges that satisfy these conditions,
shown in Figure~\ref{fig:nonpla}.
More precisely, two of the graphs are obtained from $K_{3,3}$ by adding an edge $e_1$ or $e_2$, and the other graph is $K_5$.
Thus we have the following proposition, which was mentioned in~\cite{LKLO}.

\begin{proposition} \label{prop:planar}
If $\widehat{G}_{a,b}$ is planar, then $G$ is not intrinsically knotted.
Especially, if $\widehat{G}_{a,b}$ satisfies one of the following three conditions,
then $\widehat{G}_{a,b}$ is planar, 
so $G$ is not intrinsically knotted.
\begin{itemize}
\item[(1)] $|\widehat{E}_{a,b}| \leq 8$, or
\item[(2)] $|\widehat{E}_{a,b}|=9$ and $\widehat{G}_{a,b}$ is not isomorphic to $K_{3,3}$.
\item[(3)] $|\widehat{E}_{a,b}|=10$ and $\widehat{G}_{a,b}$ is not isomorphic to $K_5$, $K_{3,3}+e_1$ and $K_{3,3}+e_2$.
\end{itemize}
\end{proposition}

\section{Restoring and twin restoring methods} \label{sec:restoring}

In this section we review the restoring method,
which we introduced in~\cite{KMO2} and will use frequently in this paper. 
We also introduce a similar technique that we call the twin restoring method.

We will find all candidate bipartite intrinsically knotted graphs with 23 edges.
To prove the main theorem, we distinguish several cases according to the combinations of degrees of 
all vertices and further sub-cases according to connections of some of the 23 edges.
Let $G$ be a bipartite graph with 23 edges with some distinct vertices $a$ and $b$.
Figure~\ref{fig:rest}(a) gives an example where $a$ and $b$ are $a_1$ and $a_2$. 
As in the figure, we assume that the degree of every vertex as well as information about 
certain edges, including all edges incident to $a$ and $b$, is known.

First, we examine the number of the edges of the graph $\widehat{G}_{a,b}$.
If it has at most eight edges, then it is planar and so $G$ cannot be intrinsically knotted by Proposition~\ref{prop:planar}.
Even if it has more edges, $G$ is rarely intrinsically knotted.
Especially if it has 9 edges, $\widehat{G}_{a,b}$ must be isomorphic to $K_{3,3}$ in order for $G$ to be intrinsically knotted.
In this case, $G_{a,b}$, being a subdivision of $K_{3,3}$, has exactly six vertices with degree 3 and, possibly, additional vertices of degree 2.
The {\em restoring method\/} is a way to find candidates for such a $G_{a,b}$ as shown in Figure~\ref{fig:rest}(b) and (c).
Finally we recover $G$ from $G_{a,b}$ by restoring the deleted vertices and edges.
$$ \stackbelow{\widehat{G}_{a,b}}{K_{3,3}}  \ \ \rightarrow \ \ G_{a,b} \ \ \rightarrow \ \ G $$

Sometimes the restoring method applied to $G_{a,b}$ for only one pair of vertices $\{a,b\}$
does not give sufficient information to construct the graph $G$.
In this case, we apply the restoring method to two graphs $G_{a,b}$ and $G_{a',b'}$ simultaneously
for different pairs of vertices.
We call this method the {\em twin restoring method\/}.

\subsection{An example of the restoring method with 9 edges} \label{subsec:restoring1} \ 

\begin{figure}[h]
\includegraphics[scale=1]{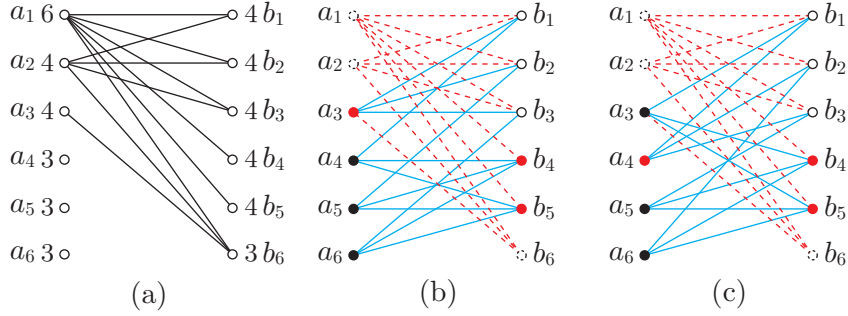}
\caption{Restoring method}
\label{fig:rest}
\end{figure}

As an example, suppose that $A$ consists of one degree 6 vertex, two degree 4 vertices and three degree 3 vertices,
and $B$ consists of five degree 4 vertices and one degree 3 vertex with edge information as shown in Figure~\ref{fig:rest}(a).
In the figure, the vertices are labeled by $a_1, \dots, a_6, b_1,\dots,b_6$ and the numbers near vertices indicate their degrees.

In this case, $G_{a_1,a_2}$ has six degree 3 vertices $a_3,a_4,a_5,a_6,b_4,b_5$ and 
three degree 2 vertices $b_1,b_2,b_3$.
Now we examine the number of the edges $|\widehat{E}_{a_1,a_2}|$ of the graph $\widehat{G}_{a_1,a_2}$.
Since $NE(a_1,a_2) = 10$, $NV_3(a_1,a_2) = 1$ and $NV_4(a_1,a_2) = 3$,
the count equation gives $|\widehat{E}_{a_1,a_2}| = 9$.

We now assume that $\widehat{G}_{a_1,a_2}$ is isomorphic to $K_{3,3}$.
As the bipartition of $K_{3,3}$, we assign the bipartition $C$ (black vertices) and 
$D$ (red vertices) for six degree 3 vertices of $G_{a_1,a_2}$.
Since all four vertices $a_3,a_4,a_5,a_6$ have degree 3,
$b_4$ is not adjacent to $b_5$ ($b_4 \nadj b_5$) in $\widehat{G}_{a_1,a_2}$.
This implies that $b_4$ and $b_5$ should be in the same partition, say $D$.
Without loss of generality, the remaining vertex of $D$ is either $a_3$ or $a_4$ 
(indeed, the three vertices $a_4, a_5, a_6$ are isomorphic).
Compare Figures~\ref{fig:rest}(b) and (c). 

In the first case, $C = \{a_4,a_5,a_6\}$.
The three edges of $\widehat{G}_{a_1,a_2}$ connecting $a_3$ to $C$ 
inevitably pass through the three degree 2 vertices $b_1, b_2, b_3$.
The three edges of $\widehat{G}_{a_1,a_2}$ incident to $b_4$ (or $b_5$)
are directly connected to $C$. 
This $G_{a_1,a_2}$ is drawn by the solid edges in the figure.
By restoring the deleted vertices and dotted edges, we recover $G$.
In the second case, $C = \{a_3,a_5,a_6\}$.
Then the three edges of $\widehat{G}_{a_1,a_2}$ connecting $a_4$ and $C$ 
passes through three degree 2 vertices $b_1, b_2, b_3$.
The remaining arguments are similar to the first case.

\subsection{An example of the restoring method with 10 edges} \label{subsec:restoring2} \ 

\begin{figure}[h]
\includegraphics[scale=1]{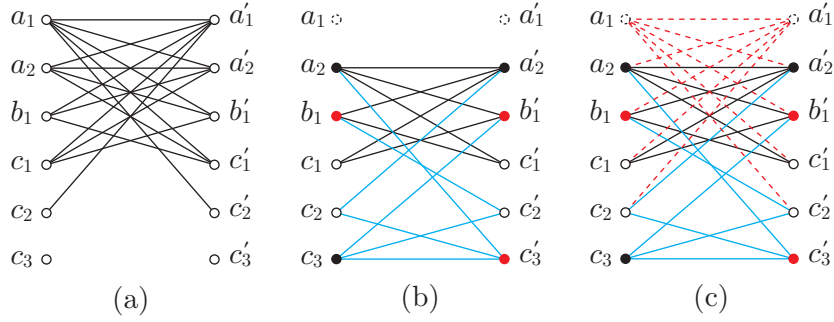}
\caption{Restoring method with 10 edges}
\label{fig:rest2}
\end{figure}

Even when $\widehat{G}_{a,b}$ has 10 edges, we can still apply the restoring method.
As an example, suppose both $A$ and $B$ consist of two degree 5 vertices, one degree 4 vertex 
and three degree 3 vertices with edge information and vertex labelling as drawn in Figure~\ref{fig:rest2}(a).
In this case, $G_{a_1,a'_1}$ has two degree 4 vertices $a_2$ and $a'_2$,
four degree 3 vertices $b_1$, $c_3$, $b'_1$ and $c'_3$,
and four degree 2 vertices $c_1$, $c_2$, $c'_1$ and $c'_2$.
Then $|\widehat{E}_{a_1,a'_1}| = 10$.

We now assume that $\widehat{G}_{a_1,a'_1}$ is isomorphic to one of $K_{3,3}+e_1$ or $K_{3,3}+e_2$.
Since $a_2$ and $a'_2$ are mutually adjacent to $b_1$ and $b'_1$ in $\widehat{G}_{a_1,a_2}$,
$a_2$ and $a'_2$ are contained in the same partition $C$ and $b_1$ and $b'_1$ are in $D$.
Therefore $\widehat{G}_{a_1,a_2}$ is isomorphic to $K_{3,3}+e_1$.
Without loss of generality, $c_3$ is contained in $C$ and $c'_3$ is contained in $D$.
Obviously $b_1$ and $c_3$ are adjacent in $\widehat{G}_{a_1,a'_1}$ passing through $c'_2$
and $a'_2$ and $c'_3$ are adjacent passing through $c_2$.
The remaining connections are drawn in Figure~\ref{fig:rest2}(b).
To recover $G$, restore the deleted vertices and dotted edges and so we get the graph $G$
as drawn in Figure~\ref{fig:rest2}(c).

\subsection{An example of the twin restoring method} \label{subsec:restoring3} \ 

\begin{figure}[h]
\includegraphics[scale=1]{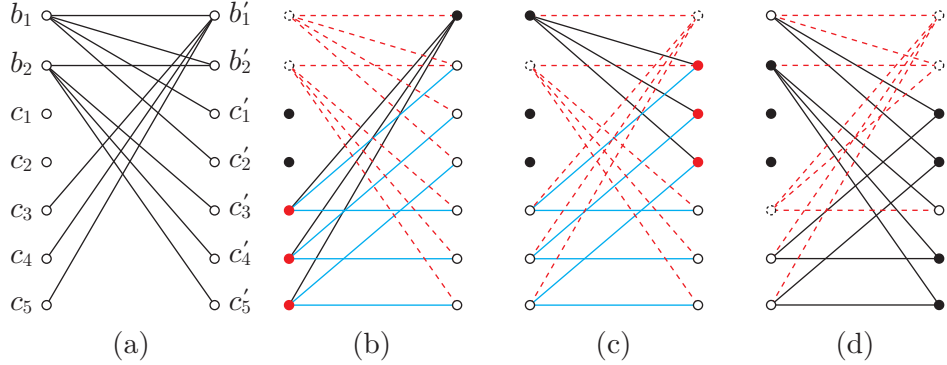}
\caption{Twin restoring method}
\label{fig:rest3}
\end{figure}

As an example, suppose that both $A$ and $B$ consist of two degree 4 vertices and five degree 3 vertices 
with vertex labelling and partial edge information as drawn in Figure~\ref{fig:rest3}(a).
In this case, we apply the restoring method to two graphs $G_{b_1,b_2}$ and $G_{b_2,b'_1}$ simultaneously.
These two graphs have the bipartitions assigned as in Figure~\ref{fig:rest3}(b) and (c).
By considering the bipartition in $G_{b_1,b_2}$, 
each of $c'_3, c'_4$ and $c'_5$ must be adjacent to exactly one of $c_3, c_4$ and $c_5$.
Furthermore, by considering the bipartition in $G_{b_2,b'_1}$, 
each of $c_3, c_4$ and $c_5$ must be adjacent to at least one of $c'_3, c'_4$ and $c'_5$.
From these two facts, we assume that $c_3 \adj c'_3$, $c_4 \adj c'_4$ and $c_5 \adj c'_5$.
Without loss of generality, we further assume that $c_3 \adj b'_2$, $c_4 \adj c'_1$, $c_5 \adj c'_2$.
In Figure~\ref{fig:rest3}(d), since $b'_2$ is adjacent to $c_1$ or $c_2$, $\widehat{G}_{b'_1,b'_2}$ has 9 edges and a 5-cycle $(b_2 c'_4 c'_1 c'_2 c'_5)$.
Since it is not isomorphic to $K_{3,3}$, $G$ is not intrinsically knotted.

\section{$G$ contains a vertex with degree 6 or more} \label{sec:g6}

Throughout this paper, we assume that $G$ is a bipartite intrinsically knotted graph with 23 edges.
In this section we assume there is a vertex $a$ in $A$ of maximal degree, with $\deg(a) \geq 6$.
We conclude that there are no size 23 bipartite intrinsically knotted graphs in this case.
Let $a'$ be a vertex in $A \setminus \{a\}$ with maximal degree.
Since $G$ has 23 edges and has vertices with degree at least 3,
$A$ and $B$ have at most seven vertices.
Therefore $\deg(a)$ is 6 or 7, and $\deg(a') \geq 4$.

Suppose $B$ has seven vertices.
Then we will say that $B$ has a 5333333 or 4433333 degree combination,
meaning either a single vertex of degree 5 or two vertices of degree 4, with the remaining vertices all 
of degree 3.
In the first case, by the count equation, 
$|\widehat{E}_{a,a'}| \leq 8$ in $\widehat{G}_{a,a'}$ since $NE(a,a') \geq 10$ and $|V_3(a)| \geq 5$.
By Proposition~\ref{prop:planar}, this contradicts $G$ being intrinsically knotted.

\begin{figure}[h]
\includegraphics[scale=1]{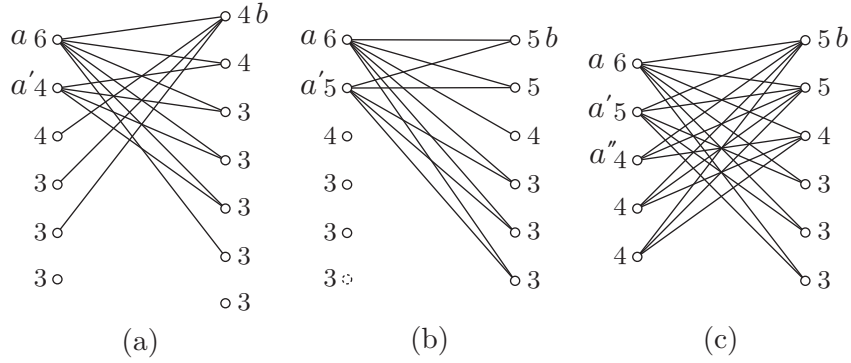}
\caption{Example of the case of $\text{deg}(a)=6$}
\label{fig:exdeg6}
\end{figure}

In the second case, 
if $\deg(a) = 7$ then $|V_3(a)|=5$, and so $|\widehat{E}_{a,a'}| \leq 7$.
If $\deg(a) = 6$ and $\deg(a') \geq 5$, then $|\widehat{E}_{a,a'}| \leq 8$.
So we can assume $\deg(a) = 6$ and $\deg(a') = 4$.
Then $A$ consists of $a$, $a'$, one more vertex with degree 4 and three vertices with degree 3.
Let $b$ be a vertex in $B$ with degree 4.
If $a' \adj b$, $|\widehat{E}_{a,a'}| \leq 8$
since $|V_3(a)|+NV_4(a,a') \geq 5$.
Otherwise, $|\widehat{E}_{a,b}| \leq 8$ since $|V_3(a)|+|V_3(b)| \geq 6$.
See Figure~\ref{fig:exdeg6}(a) for an example.

Now we assume that $B$ has six vertices, and so $a$ has degree 6.
If $a'$ has degree 6, then $NE(a,a') = 12$ and $NV_3(a,a') + NV_4(a,a') \geq 3$.
So $|\widehat{E}_{a,a'}| \leq 8$.

Suppose $\deg(a') = 5$.
If $B$ has either at least four vertices with degree 3 or at least five vertices with degree 3 or 4,
then $NV_3(a,a') + NV_4(a,a') \geq 4$, and so $|\widehat{E}_{a,a'}| \leq 8$.
So we can assume, $B$ has both at most three vertices with degree 3 and at most four vertices with degree 3 or 4,
meaning $B$ must have 554333 degree combination.
If $a'$ is adjacent to a degree 4 vertex in $B$, 
$NV_3(a,a') + NV_4(a,a') = 4$, and so $|\widehat{E}_{a,a'}| \leq 8$.
Therefore, we assume that $a'$ is not adjacent to a degree 4 vertex in $B$ as in Figure~\ref{fig:exdeg6}(b).
Here $b$ is a degree 5 vertex in $B$.
Now, $A$'s degree combination is one of 65543, 65444 or 653333.
If $A$ has a 65543 or 653333 degree combination,
then $|\widehat{E}_{a,b}| \leq 9$ with a degree 4 vertex $a'$ in $\widehat{G}_{a,b}$.
Since it is not isomorphic to $K_{3,3}$, $G$ is not intrinsically knotted.
If $A$ has a 65444 degree combination,
then the two degree 5 vertices in $B$ are adjacent to all vertices in $A$
and the degree 4 vertex in $B$ is adjacent to all vertices except $a'$ as in Figure~\ref{fig:exdeg6}(c).
For a degree 4 vertex $a''$ in $A$,
$\widehat{G}_{a,a''}$ has at most 9 edges and a degree 4 vertex $a'$.

The remaining case is that $\deg(a') = 4$, and so $A$ has a 644333 degree combination.
Let $b$ be a vertex in $B$ with maximal degree.
If $\deg(b) = 6$, then $|\widehat{E}_{a,b}| \leq 8$.
If $\deg(b) = 5$  and $B$ has at least three degree 3 vertices, 
then $NE(a,b) = 10$ and $NV_3(a,b) \geq 5$, implying $|\widehat{E}_{a,b}| \leq 8$.
Now consider the case that either $\deg(b) = 4$ or 
$\deg(b) = 5$ along with the condition that $B$ has at most two degree 3 vertices.
Then $B$ has a 444443 or 544433 degree combination.

Assume that $A$ and $B$ have 644333 and 444443 degree combinations, respectively.
If $a'$ is not adjacent to the unique degree 3 vertex in $B$, then $|\widehat{E}_{a,a'}| \leq 8$.
So the two degree 4 vertices in $A$ are adjacent to the degree 3 vertex in $B$, 
and without loss of generality we have the connecting combination 
as in Figure~\ref{fig:rest}(a).
The rest of process follows the restoring method, 
discussed in Subsection~\ref{subsec:restoring1} as an example.
Eventually we obtain the two graphs for $G$ shown in Figures~\ref{fig:rest}(b) and (c).
In Figure~\ref{fig:rest}(b), $\widehat{G}_{a_1,a_4}$ is planar.
In Figure~\ref{fig:rest}(c), $\widehat{G}_{a_1,a_3}$ has at most 9 edges and a 2-cycle $(a_5 a_6)$.

\begin{figure}[h]
\includegraphics[scale=1]{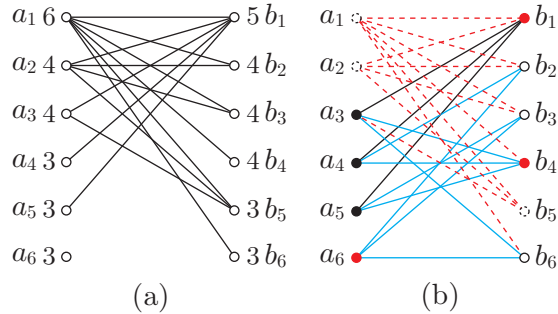}
\caption{The case of $\text{deg}(a_1)=6$ and the restoring method}
\label{fig:deg6}
\end{figure}

Now consider the final case where $B$ has a 544433 degree combination.
We label the vertices in descending order of their vertex degree as in Figure~\ref{fig:deg6}(a).
If $|V_3(b_1)| =3$ then $|\widehat{E}_{a_1,b_1}| \leq 8$.
Without loss of generality, we may assume that $b_1 \nadj a_6$.
If $|V_3(a_2)| =0$ (similarly for $a_3$) then $|\widehat{E}_{a_1,a_2}| \leq 8$.
Also if $|V_3(a_2)| =|V_3(a_3)|=2$ 
then $\widehat{G}_{a_1,b_1}$ has at most 9 edges and a 2-cycle $(a_2 a_3)$.
Therefore one of them, say $a_2$, is adjacent to exactly one degree 3 vertex in $B$. 
Without loss of generality, $a_2$ is adjacent to $b_1,b_2,b_3,b_5$.
Now $a_3 \adj b_5$, for otherwise, $NV_3(a_1,a_2) + NV_4(a_1,a_2) = 4$ 
and $NV_Y(a_1,a_2) = 1$, and so $|\widehat{E}_{a_1,a_2}| \leq 8$.
Here $V_Y(a_1,a_2)$ includes the degree 3 vertex in $A$, which is adjacent to $b_5$. 

If $a_3$ is adjacent to $b_2$ (or $b_3$), then $|\widehat{E}_{a_1,b_1}| \leq 9$ 
with a 3-cycle $(a_2 a_3 b_2)$.
So $a_3$ is adjacent to $b_4$ and $b_6$.
By applying the restoring method, we construct $G_{a_1,a_2}$ as drawn in Figure~\ref{fig:deg6}(b).
Finally we recover $G$ by restoring the deleted vertices and dotted edges.
For the graph $G$, $\widehat{G}_{a_1,a_3}$ has 9 edges and a 3-cycle $(a_2 b_2 b_3)$,
showing that $G$ is not intrinsically knotted.

\section{Both $A$ and $B$ contain degree 5 vertices} \label{sec:ab5}

In this section we assume $G$ has maximal degree 5 and both $A$ and $B$ have degree 5 vertices.
We find two graphs for Theorem~\ref{thm:main} in Subsection~\ref{subsec:110} (see Figure~\ref{fig:b320}).
Both are formed by adding an edge to Cousin 110 of $E_9+e$ family.
Let $A_n$ denote the set of vertices in $A$ with degree $n = 3,4,5$ and $[A] = [|A_5|, |A_4|,|A_3|]$.
The possible cases for $[A]$ are $[4,0,1]$, $[3,2,0]$, $[2,1,3]$, $[1,3,2]$ and $[1,0,6]$.
Similarly, define $B_n$ and $[B]$.
Without loss of generality, we may assume that $|A_5| \geq |B_5|$,
and furthermore, if $|A_5| = |B_5|$ then $|A_4| \geq |B_4|$.

We distinguish fifteen cases of all possible combinations of $[A]$ and $[B]$, which we 
treat in the following seven subsections.
To simplify the notation, vertices in $A_5$, $A_4$, $A_3$, $B_5$, $B_4$ and $B_3$ 
are denoted by $\{a_i\}$, $\{b_i\}$, $\{c_i\}$, $\{a'_i\}$, $\{b'_i\}$ and $\{c'_i\}$, respectively.

\subsection{$[A]=[4,0,1]$ or $[3,2,0]$, and $[B]=[4,0,1]$ or $[3,2,0]$} \label{subsec:110} \ 

If both are $[4,0,1]$, the four degree 5 vertices in $A$ must all be adjacent to the unique degree 3 vertex in $B$,
which is impossible.

\begin{figure}[h]
\includegraphics[scale=1]{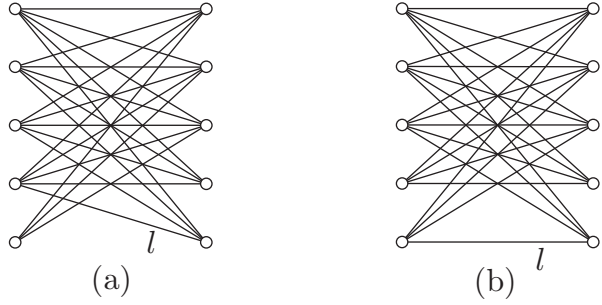}
\caption{Two bipartite intrinsically knotted graphs}
\label{fig:b320}
\end{figure}

Suppose instead, $[A]=[4,0,1]$ or $[3,2,0]$, and $[B]=[3,2,0]$.
Both cases are uniquely realized as the two graphs in Figure~\ref{fig:b320}, 
which are obtained from Cousin 110 of $E_9+e$ family by adding an edge $l$.
Cousin 110 of $E_9+e$ is intrinsically knotted~\cite{GMN}.
These are the first two bipartite intrinsically knotted graphs of Theorem~\ref{thm:main}.

\subsection{$[A]=[4,0,1]$ or $[3,2,0]$, and $[B]=[2,1,3]$} \ 

\begin{figure}[h]
\includegraphics[scale=1]{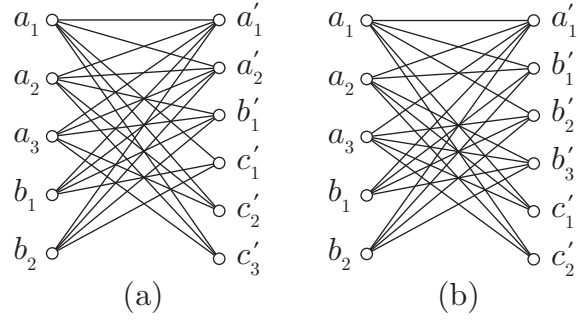}
\caption{Case of $[A]=[3,2,0]$, and $[B]=[2,1,3]$ or $[1,3,2]$}
\label{fig:a51}
\end{figure}

First assume that a degree 3 vertex $c'_1$ in $B$ is adjacent to at most one degree 5 vertex in $A$.
In this case, $[A]=[3,2,0]$ and $c'_1$ must be adjacent to one degree 5 vertex $a_1$ and 
two degree 4 vertices in $A$.
Furthermore $a_2$ and $a_3$ are adjacent to all vertices in $B$ except $c'_1$.
If $a_1 \adj b'_1$, then $|\widehat{E}_{a_1,a_2}| \leq 9$ and $a_3$ has degree larger than 3 
in $\widehat{G}_{a_1,a_2}$.
If $a_1 \nadj b'_1$, then $|\widehat{E}_{a_2,b_1}| \leq 10$ and $a_1$ has degree 5 
in $\widehat{G}_{a_2,b_1}$ as drawn in Figure~\ref{fig:a51}(a).

So every degree 3 vertex in $B$ is adjacent to at least two degree 5 vertices in $A$.
Assume that $c'_1 \adj \{a_1, a_2\}$, and furthermore $b'_1 \adj a_1$.
If $a_3 \nadj c'_1$, then $\widehat{G}_{a_1,a_3}$ has at most 9 edges and 
a vertex $a_2$ with degree larger than 3.
So every degree 5 vertex in $A$ is adjacent to all degree 3 vertices in $B$, and so $[A]=[3,2,0]$.
In this case no degree 4 vertex in $A$ can be adjacent to all degree 3 vertices in $B$,
but no such graph $G$ is possible.

\subsection{$[A]=[4,0,1]$ or $[3,2,0]$, and $[B]=[1,3,2]$} \ 

If a degree 5 vertex $a_1$ in $A$ is adjacent to all degree 4 vertices in $B$,
then $|\widehat{E}_{a_1,a_2}| \leq 9$ and $a_3$ has degree larger than 3 in $\widehat{G}_{a_1,a_2}$.
Suppose instead each degree 5 vertex in $A$ is adjacent to two degree 3 vertices 
and two among the three degree 4 vertices in $B$.
So $[A]=[3,2,0]$ and $G$ is uniquely realized as in Figure~\ref{fig:a51}(b).
Since $\widehat{G}_{a_1,b_2}$ has 10 edges with a degree 5 vertex, it is planar.

\subsection{$[A]=[2,1,3]$ and $[B]=[2,1,3]$.}\ 

If a degree 5 vertex $a_1$ is adjacent to all three degree 3 vertices in $B$,
then $|\widehat{E}_{a_1,a'_1}| \leq 9$ and
$a_2$ has degree larger than 3 in $\widehat{G}_{a_1,a'_1}$.
The same argument applies to vertices $a_2, a'_1$ and $a'_2$.
Therefore there are vertices $c_1$ and $c'_1$ adjacent to both degree 5 vertices on the other side
as in Figure~\ref{fig:rest2}(a).
As in the figure, we also have the condition $b_1 \adj c'_1$ (similarly $b'_1 \adj c_1$).
For, if $a_2 \adj c'_3$ and $b_1 \nadj c'_1$, $|\widehat{E}_{a_1,a_2}| \leq 8$
because $V_Y(a_1,a_2)$ is not empty.
Or, if $a_2 \adj c'_2$ and $b_1$ is not adjacent to both $c'_1$ and $c'_2$, 
then $b_1 \adj \{ a'_1,a'_2,b'_1,c'_3 \}$.
This implies that $\widehat{G}_{a_1,b_1}$ has 10 edges and a degree 5 vertex $a_2$.

The rest of process follows the restoring method, 
as described in Subsection~\ref{subsec:restoring2} as an example.
This leads to the graph $G$ as drawn in Figure~\ref{fig:rest2}(c).
Then $\widehat{G}_{a'_1,a'_2}$ has at most 9 edges and a 3-cycle $(a_2 c_3 b'_1)$.

\subsection{$[A]=[2,1,3]$ and $[B]=[1,3,2]$.}\ 

If there is a degree 5 vertex $a_1$ which is not adjacent to $a'_1$,
$NV_3(a_1,a'_1)=5$, and so $|\widehat{E}_{a_1,a'_1}| =8$.
Therefore $a'_1 \adj \{a_1, a_2\}$.

\begin{figure}[h]
\includegraphics[scale=1]{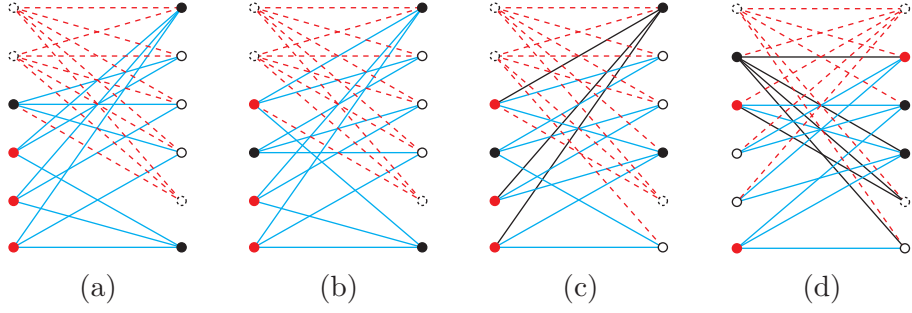}
\caption{Case of $[A]=[2,1,3]$ and $[B]=[1,3,2]$}
\label{fig:213132}
\end{figure}

First assume that $NV_4(a_1,a_2)=3$.
Let $a_1 \adj c'_1$.
To avoid $|\widehat{E}_{a_1,a_2}| \leq 8$, $NV_3(a_1,a_2)+NV_Y(a_1,a_2) \leq 1$.
This implies that $a_2 \adj c'_1$ and $V_Y(a_1,a_2)$ is empty (so $b_1 \adj c'_1$).
Now we apply the restoring method.
According to the bipartition choice of $\widehat{G}_{a_1,a_2}$, 
we construct two graphs for $G_{a_1,a_2}$ as drawn in Figure~\ref{fig:213132}(a) and (b).
Finally we recover $G$ by restoring the deleted vertices and dotted edges.
In Figure~\ref{fig:213132}(a) and (b), 
we find planar graphs $\widehat{G}_{a_1,c_1}$ and $\widehat{G}_{a_1,b_1}$, respectively.

Assume that $NV_4(a_1,a_2)=2$, which are $b'_1,b'_2$.
Obviously, $NV_3(a_1,a_2)=2$.
If both $a_1$ and $a_2$ are not adjacent to $b'_3$, 
then $NV_Y(a_1,a_2) \geq 1$, implying $|\widehat{E}_{a_1,a_2}| \leq 8$.
Now we assume that $a_1 \adj b'_3$ and $a_2 \nadj b'_3$,
and assume further that $a_1 \adj c'_1$ and $a_2 \adj \{c'_1, c'_2\}$.
If $c'_1 \adj c_i$, then $NV_Y(a_1,a_2)=1$, and so $c'_1 \adj b_1$.
If $a'_1 \nadj b_1$, then $\widehat{G}_{a_2,a'_1}$ has 9 edges and a degree 4 vertex $a_1$, 
and so we assume $a'_1 \adj \{b_1,c_1,c_2\}$ as in Figure~\ref{fig:213132}(c).
Using the restoring method, we construct $G_{a_1,a_2}$, and then recover $G$.
In this graph, $\widehat{G}_{a_1,b_1}$ is planar.

Finally we have $NV_4(a_1,a_2)=1$, which is $b'_1$.
Then we may assume that $a_1 \adj \{a'_1,b'_1,b'_2,c'_1,c'_2\}$ and $a_2 \adj \{a'_1,b'_1,b'_3,c'_1,c'_2\}$.
If $b_1 \adj \{a'_1,b'_1,b'_2,b'_3\}$, 
then $\widehat{G}_{a_1,b_1}$ has 10 edges and a degree 5 vertex $a_2$.
Therefore, assume $b_1 \adj c'_1$.
If $a'_1 \nadj b_1$, then $\widehat{G}_{a_1,a'_1}$ has 9 edges and a degree 4 vertex $a_2$.
Thus we assume $a'_1 \adj \{b_1,c_1,c_2\}$.
Since $\widehat{G}_{a_1,a'_1}$ has 10 edges and exactly two degree 4 vertices $a_2$ and $b'_3$,
it is isomorphic to either $K_{3,3}+e_1$ or $K_{3,3}+e_2$.
However, $K_{3,3}+e_2$ is not possible for $\widehat{G}_{a_1,a'_1}$.
Therefore $\widehat{G}_{a_1,a'_1}$ is isomorphic to $K_{3,3}+e_1$.
Using the restoring method, 
we construct $G_{a_1,a'_1}$ as drawn in Figure~\ref{fig:213132}(d), and then recover $G$.
Since $NV_3(a_1,a_2)+NV_4(a_1,a_2)+NV_Y(a_1,a_2)=4$, 
$\widehat{G}_{a_1,a_2}$ has 9 edges and a 3-cycle $(a'_1 c_1 c_2)$.

\subsection{$[A]=[1,3,2]$ and $[B]=[1,3,2]$.}\ 

\begin{figure}[h!]
\includegraphics{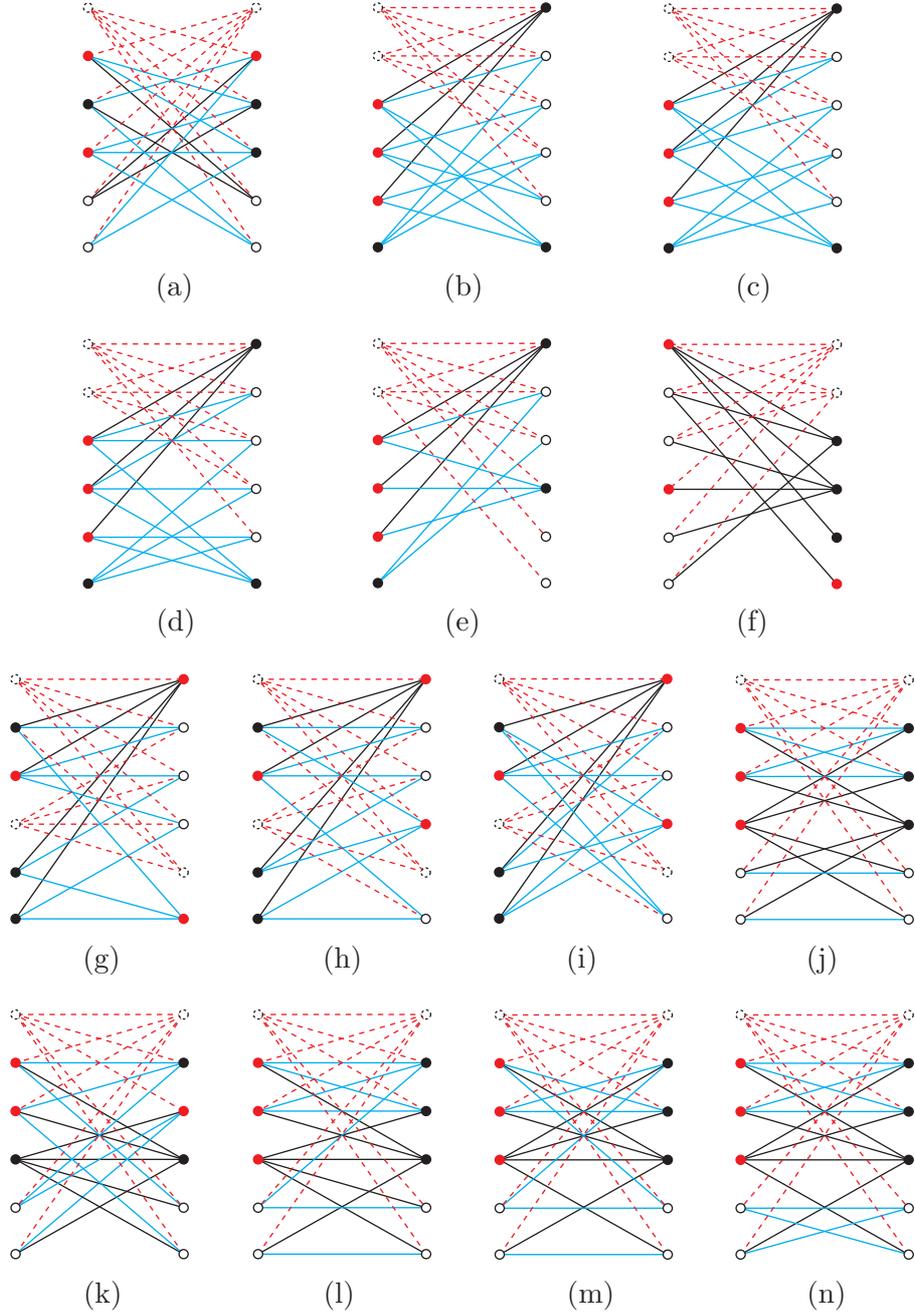}
\caption{Case of $[A]=[1,3,2]$ and $[B]=[1,3,2]$.}
\label{fig:132132}
\end{figure}

First assume that $a_1 \nadj a'_1$.
If there is a degree 4 vertex $b_1$ in $A$ which is not adjacent to degree 3 vertices in $B$,
then $NV_4(a_1,b_1)=3$, implying that $\widehat{G}_{a_1,b_1}$ has 9 edges and a degree 4 vertex $a'_1$.
Therefore we may assume that $\{b_1,b_2\} \adj c'_1$, and similarly $\{b'_1,b'_2\} \adj c_1$.
Using the restoring method, we construct $G_{a_1,a'_1}$ as in Figure~\ref{fig:132132}(a).
After recovering $G$, we have a planar graph $\widehat{G}_{b_1,b_2}$.

Now assume that $a_1 \adj a'_1$.
We distinguish into three cases.
The first case is that $a_1 \adj \{b'_1,b'_2,b'_3,c'_1\}$ and $a'_1 \adj \{b_1,b_2,b_3,c_1\}$.
There is a vertex, say $b_1$, among the degree 4 vertices in $A$ such that $b_1 \nadj c'_1$.
We may assume either $b_1 \adj \{b'_1,b'_2,b'_3\}$ or $b_1 \adj \{b'_1,b'_2,c'_2\}$.
In the former case, 
using the restoring method, we construct three graphs $G_{a_1,b_1}$ as drawn in 
Figure~\ref{fig:132132}(b), (c) and (d).
In all cases, after recovering $G$, we have planar graphs $\widehat{G}_{b_1,b_2}$.
In the latter case, 
we partially construct $\widehat{G}_{a_1,b_1}$ as in Figure~\ref{fig:132132}(e).
In the figure, the bipartition is determined by the connection of $a'_1$ so that
$b'_3 \adj \{b_2,b_3,c_1\}$.
To be $K_{3,3}+e_1$, the two degree 4 vertices $b_2, b_3$ must be connected through a degree 2 vertex.
If this vertex is $b'_1$ (or similarly $b'_2$), then $b'_1 \adj \{b_1,b_2,b_3\}$,
which is the same as the former case with exchanging the vertices in $A$ and $B$.
This implies that $c_2 \adj \{b'_1,b'_2\}$ and further we may assume $b_2 \adj b'_1$.
Now we instead consider $\widehat{G}_{a'_1,b'_1}$ as in Figure~\ref{fig:132132}(f).
Then the bipartition is determined by the connection of $a_1$,
but the two degree 3 vertices $c'_1, c'_2$ cannot be connected through any degree 2 vertex.

The second case is that $a_1 \adj \{b'_1,b'_2,b'_3,c'_1\}$ and $a'_1 \adj \{b_1,b_2,c_1,c_2\}$.
If $b_3 \nadj c'_1$, then $b_3 \adj \{b'_1,b'_2,b'_3,c'_2\}$,
implying that $\widehat{G}_{a_1,b_3}$ has 9 edges and a degree 4 vertex $a'_1$.
Thus we have $b_3 \adj c'_1$.
If $c'_1$ is adjacent to a degree 3 vertex in $A$, then $V_Y(a_1,b_3)$ is not empty,
implying that $\widehat{G}_{a_1,b_3}$ has 9 edges and a degree 4 vertex $a'_1$.
Therefore, $b_1 \adj c'_1$.
Now in either case of $b_3 \adj \{b'_1,b'_2,b'_3\}$ or $b_3 \adj \{b'_1,b'_2,c'_2\}$,
using the restoring method, we have three graphs $\widehat{G}_{a_1,b_3}$ 
as in Figure~\ref{fig:132132}(g), (h) and (i).
In all cases, after recovering $G$, we have planar graphs $\widehat{G}_{a_1,b_1}$.

Finally we consider the third case, where $a_1 \adj \{b'_1,b'_2,c'_1,c'_2\}$ and $a'_1 \adj \{b_1,b_2,c_1,c_2\}$.
In this case, we construct $\widehat{G}_{a_1,a'_1}$ which has 10 edges with exactly two degree 4 vertices
$b_3, b'_3$, so it must be isomorphic to either $K_{3,3}+e_1$ or $K_{3,3}+e_2$.
Then we have four cases as follows;
(1) $b_3 \nadj b'_3$, 
(2) $b_3 \nadj b'_1$ and $b'_3 \nadj b_1$, 
(3) $b_3 \nadj b'_1$ and $b'_3 \nadj c_1$, and
(4) $b_3 \nadj c'_1$ and $b'_3 \nadj c_1$.
In the figure, (j) indicates the case (1), (k) and (l) indicate the case (3), 
(m) and (n) indicate the case (4), and no graph satisfying case (2) can be constructed.
We find a planar graph $\widehat{G}_{a_1,b_3}$ in (j), and planar graphs $\widehat{G}_{a_1,b_1}$
for the remaining cases.

\subsection{$[A]$ is one of the five cases, and $[B]=[1,0,6]$} \ 

First assume that $[A]=[4,0,1]$ or $[3,2,0]$.
If $NV_3(a_i,a_j) \geq 5$ for some $i$ and $j$, then $|\widehat{E}_{a_i,a_j}| \leq 8$.
Suppose instead, three degree 5 vertices in $A$ are adjacent to the unique degree 5 vertex and 
the same four degree 3 vertices in $B$.
It is not possible to construct such a graph $G$.

If $[A]=[2,1,3]$ or $[1,0,6]$, $NV_3(a_1,a'_1) \geq 6$, implying $|\widehat{E}_{a_1,a'_1}| \leq 8$.

Finally, assume that $[A]=[1,3,2]$. 
If $a_1 \nadj a'_1$, $NV_3(a_1,a'_1) = 7$, implying $|\widehat{E}_{a_1,a_2}| \leq 6$.
Now we assume that $a_1 \adj \{a'_1, c'_1, c'_2, c'_3, c'_4\}$.
If $b_1 \adj \{c'_5,c'_6\}$ then $NV_3(a_1,b_1) = 6$, implying $|\widehat{E}_{a_1,a_2}| \leq 8$.
So none of the degree 4 vertices $b_1, b_2, b_3$ is adjacent to both $c'_5$ and $c'_6$.
This means that there are at least three edges connecting $\{c_1,c_2\}$ and $\{c'_5,c'_6\}$.
Therefore, one of them, say $c_1$, is adjacent to both $c'_5$ and $c'_6$.
Since $NV_3(a_1,c_1) = 6$, $\widehat{G}_{a_1,c_1}$ has 9 edges 
and a degree 4 vertex among $b_1, b_2$ and $b_3$.

\section{Only $A$ contains degree 5 vertices} \label{sec:a5}

In this section we assume that only $A$ contains degree 5 vertices and 
$B$ contains vertices with degree at most 4.
In Subsection~\ref{subsec:106025} we find a bipartite intrinsically knotted graph
formed by adding two edges to the Heawood graph (see Figure~\ref{fig:106025}(d)).
The possible cases for $[A]$ are $[4,0,1]$, $[3,2,0]$, $[2,1,3]$, $[1,3,2]$ and $[1,0,6]$,
and for $[B]$, $[0,5,1]$ and $[0,2,5]$.

We distinguish ten cases of possible combinations of $[A]$ and $[B]$
in the following six subsections.

\subsection{$[A]=[4,0,1]$, $[3,2,0]$ or $[2,1,3]$, and $[B]=[0,5,1]$} \ 

If two degree 5 vertices $a_1$ and $a_2$ in $A$ satisfy $V_4(a_1,a_2) \geq 4$, 
then $NV_3(a_1,a_2) + NV_4(a_1,a_2)=5$, implying $|\widehat{E}_{a_1,a_2}|=8$.
Suppose instead that $a_1 \adj \{b'_1,b'_2,b'_3,b'_4,c'_1\}$ and $a_2 \adj \{b'_1,b'_2,b'_3,b'_5,c'_1\}$.

If $[A]=[4,0,1]$ or $[3,2,0]$, 
$\widehat{G}_{a_1,a_2}$ has at most 9 edges and a vertex $a_3$ with degree larger than 3.

\begin{figure}[h!]
\includegraphics{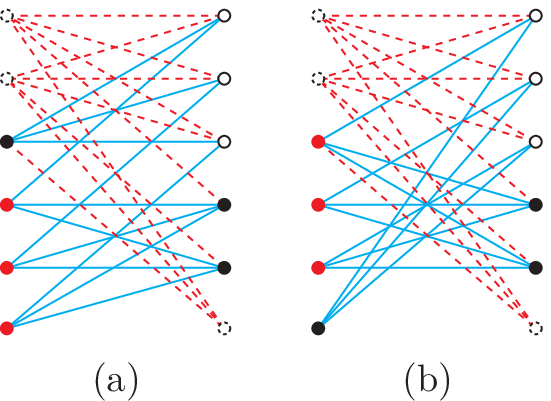}
\caption{Case of $[A]=[2,1,3]$ and $[B]=[0,5,1]$.}
\label{fig:213051}
\end{figure}

On the other hand, if $[A]=[2,1,3]$, then we have $b_1 \adj c'_1$
because, if not, $NV_Y(a_1,a_2)=1$, implying $|\widehat{E}_{a_1,a_2}|=8$.
Using the restoring method, we construct two graphs of $G_{a_1,a_2}$ 
as in Figure~\ref{fig:213051}(a) and (b).
By recovering $G$, we find that $\widehat{G}_{a_1,c_1}$ is planar in both cases.

\subsection{$[A]=[4,0,1]$, $[3,2,0]$ or $[2,1,3]$, and $[B]=[0,2,5]$} \ 

If a degree 5 vertex $a_1$ in $A$ is adjacent to either both or neither of the two degree 4 vertices in $B$,
then $NV_3(a_1,a_2) + NV_4(a_1,a_2) \geq 5$, implying $|\widehat{E}_{a_1,a_2}| \leq 8$.
Therefore each degree 5 vertex in $A$ must be adjacent to exactly one degree 4 vertex.

If $[A]=[4,0,1]$ or $[3,2,0]$, let $a_1$ and $a_2$ be two degree 5 vertices in $A$ 
that are adjacent to the same degree 4 vertex in $B$, implying $|\widehat{E}_{a_1,a_2}| \leq 8$.

If $[A]=[2,1,3]$, then we may say that $a_1 \nadj b'_1$.
In this case $NV_3(a_1,b'_1) \geq 6$, implying $|\widehat{E}_{a_1,b'_1}| \leq 8$.

\subsection{$[A]=[1,3,2]$ and $[B]=[0,5,1]$.}\ 

\begin{figure}[h!]
\includegraphics{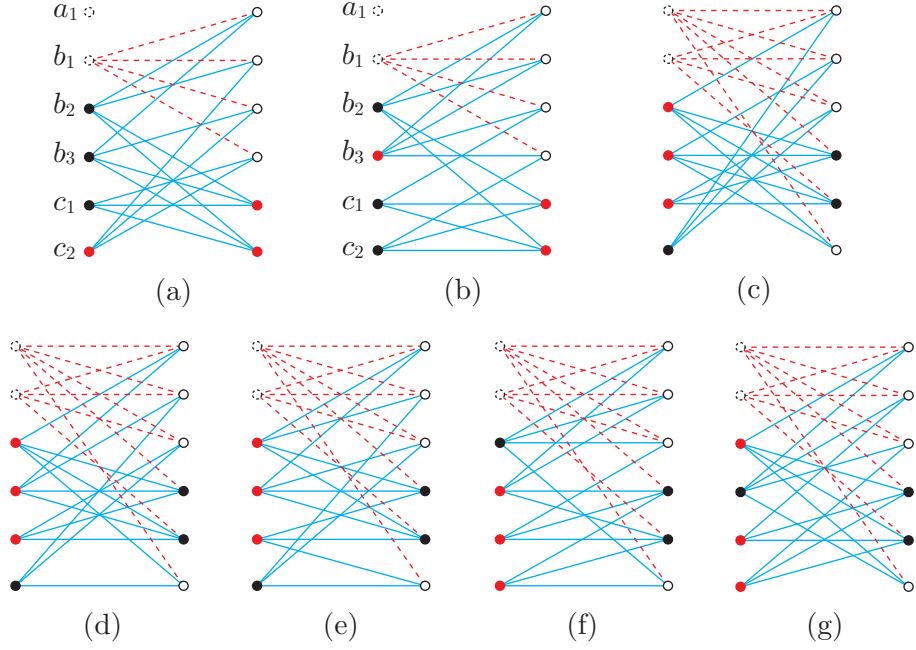}
\caption{Case of $[A]=[1,3,2]$ and $[B]=[0,5,1]$.}
\label{fig:132051}
\end{figure}

First assume that $a_1$ is adjacent to all degree 4 vertices of $B$.
Since $\widehat{G}_{a_1,b_1}$ has 10 edges and exactly two degree 4 vertices $b_2$ and $b_3$,
it is isomorphic to either $K_{3,3}+e_1$ or $K_{3,3}+e_2$.
For the first case as in Figure~\ref{fig:132051}(a), 
we recover $G \setminus \{ a_1 \}$ instead of $G$ by restoring only $b_1$ and the related edges.
In this graph, $\widehat{G}_{a_1,b_3}$ is planar.
For the second case as in Figure~\ref{fig:132051}(b), 
we similarly recover $G \setminus \{ a_1 \}$, and then $\widehat{G}_{a_1,b_2}$ is planar.

Now we assume that $a_1 \adj \{b'_1,b'_2,b'_3,b'_4,c '_1\}$.
Then we further assume that $c'_1 \nadj b_1$.
If $NV_4(a_1,b_1)=4$, 
then $\widehat{G}_{a_1,b_1}$ has at most 9 edges and a vertex $b_2$ with degree larger than 3.
So we may say that $b_1 \adj \{b'_1,b'_2,b'_3,b'_5\}$.
Since $\widehat{G}_{a_1,b_1}$ has 10 edges and exactly two degree 4 vertices $b_2$ and $b_3$,
it is isomorphic to either $K_{3,3}+e_1$ or $K_{3,3}+e_2$.
Using the restoring method, 
we construct $G_{a_1,b_1}$ as drawn in Figure~\ref{fig:132051}(c), (d), (e), (f), and (g), 
in which the first three figures correspond to $K_{3,3}+e_1$ and the remaining two figures to $K_{3,3}+e_2$.
After recovering $G$, we have planar graphs $\widehat{G}_{a_1,b_2}$ for all five cases.

\subsection{$[A]=[1,3,2]$ and $[B]=[0,2,5]$.}\ 

If $a_1$ is adjacent to all five degree 3 vertices in $B$,
$NV_3(a_1,b'_1) \geq 6$, implying $|\widehat{E}_{a_1,b'_1}| \leq 8$.

Now assume that $a_1 \adj b'_1$ and $a_1 \nadj b'_2$.
If $b'_2 \adj \{c_1, c_2\}$, then $NV_3(a_1,b'_2) \geq 6$, implying $|\widehat{E}_{a_1,b'_2}| \leq 8$.
Therefore $b'_2$ is adjacent to all three degree 4 vertices in $A$.
Similarly if $b'_1 \adj \{c_1, c_2\}$, then $\widehat{G}_{a_1,b'_1}$ has 9 edges and a degree 4 vertex $b'_2$, 
and so we assume $b'_1 \adj \{b_1,b_2\}$.
To avoid $|\widehat{E}_{a_1,b_1}| \leq 8$, $NV_3(a_1,b_1)+NV_Y(a_1,b_1) \leq 4$ because $NV_4(a_1,b_1)=1$.
So we may assume that $b_1 \adj \{c'_1,c'_2\}$, and then $c'_1 \adj b_2$ and $c'_2 \adj b_3$.
Therefore $\widehat{G}_{a_1,b'_2}$ has 9 edges and a 3-cycle $(b_1 b_2 b'_1)$.

Finally we assume that $a_1 \adj \{b'_1,b'_2,c'_1,c'_2,c'_3\}$.
If a degree 4 vertex $b_i$ in $A$ is adjacent to at most one among $\{c'_1,c'_2,c'_3\}$,
then $|\widehat{E}_{a_1,b_i}| \leq 8$.
Therefore $NV_3(a_1,b_i) \geq 2$,
and we may assume that $b_1 \adj \{c'_1,c'_2\}$, $b_2 \adj \{c'_1,c'_3\}$ and $b_3 \adj \{c'_2,c'_3\}$.
If a degree 4 vertex $b'_i$ in $B$ is adjacent to at least two degree 4 vertices in $A$, say $b_1$ and $b_2$,
then $\widehat{G}_{a_1,b_3}$ has 9 edges and a 3-cycle $(b_1 b_2 b'_1)$.
Therefore both $b'_1$ and $b'_2$ are adjacent to $c_1$ and $c_2$.
We further assume that $b'_1 \adj b_1$, 
and so $\widehat{G}_{a_1,b_1}$ has 9 edges and a 3-cycle $(c_1 c_2 b'_2)$.

\subsection{$[A]=[1,0,6]$ and $[B]=[0,5,1]$.}\ 

\begin{figure}[h!]
\includegraphics{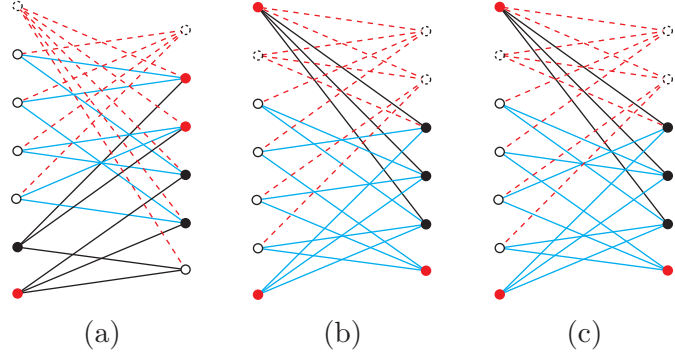}
\caption{Case of $[A]=[1,0,6]$ and $[B]=[0,5,1]$.}
\label{fig:106051}
\end{figure}

If there is a degree 4 vertex $b'_1$ in $B$ which is not adjacent to $a_1$,
then we assume that $b'_1 \adj \{c_1, c_2, c_3, c_4\}$.
If both $c_5$ and $c_6$ are adjacent to the same degree 4 vertex, say $b'_2$,
then $\widehat{G}_{b'_1,b'_2}$ has 9 edges and a degree 4 vertex $a_1$.
Therefore we assume that $c_5 \adj \{b'_2, b'_3, c'_1\}$ and $c_6 \adj \{b'_4, b'_5, c'_1\}$.
Now use the restoring method to construct $G_{a_1,b'_1}$ as drawn in Figure~\ref{fig:106051}(a). 
After recovering $G$, we have a planar graph $\widehat{G}_{a_1,c_6}$.

Suppose instead that $a_1$ is adjacent to all degree 4 vertices in $B$.
If there is a pair of degree 4 vertices, say $b'_1, b'_2$, so that $V_3(b'_1, b'_2)$ is empty,
then $\widehat{G}_{b'_1,b'_2}$ has 9 edges and a degree 4 vertex $b'_3$.
Furthermore it is impossible that all pairs of degree 4 vertices, for example $b'_1, b'_2$, 
are such that $V_3(b'_1, b'_2)$ has at least two degree 3 vertices in $A$.
Therefore we may assume that $b'_1 \adj \{c_1, c_2, c_3\}$ and $b'_2 \adj \{c_1, c_4, c_5\}$.
Using the restoring method,
we construct $G_{b'_1,b'_2}$ as drawn in Figure~\ref{fig:106051}(b) and (c). 
After recovering $G$, we have planar graphs $\widehat{G}_{b'_1,b'_3}$ in both cases.

\subsection{$[A]=[1,0,6]$ and $[B]=[0,2,5]$.}\label{subsec:106025} \ 

If there is a degree 4 vertex $b'_1$ in $B$ which is not adjacent to $a_1$,
then $NV_3(a_1,b'_1) \geq 7$, implying $|\widehat{E}_{a_1,b'_1}| \leq 7$.
Therefore we assume that $a_1 \adj \{b'_1, b'_2, c'_1, c'_2, c'_3\}$.

\begin{figure}[h!]
\includegraphics{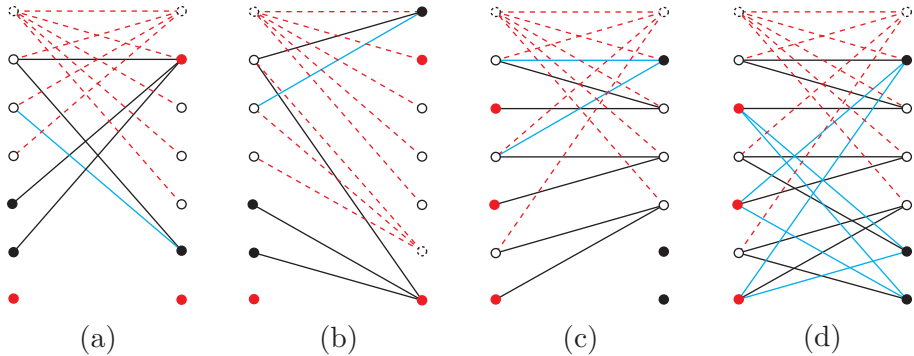}
\caption{Case of $[A]=[1,0,6]$ and $[B]=[0,2,5]$.}
\label{fig:106025}
\end{figure}

Suppose that $c_1$ is not adjacent to $c'_1, c'_2, c'_3$.
Then we distinguish two cases: $c_1 \adj \{b'_1, b'_2, c'_4\}$ or $c_1 \adj \{b'_1, c'_4, c'_5\}$.
In the first case, if both $b'_1, b'_2$ are adjacent to a degree 3 vertex $c_2$,
then $\widehat{G}_{a_1,c'_4}$ has at most 9 edges and a 3-cycle $(c_2 b'_1 b'_2)$.
So we may assume that $b'_1 \adj \{c_2, c_3\}$ and $b'_2 \adj \{c_4, c_5\}$.
Now we partially construct $\widehat{G}_{a_1,b_1}$ as in Figure~\ref{fig:106025}(a).
In the figure, the bipartition is determined by the connection of $b'_2$, 
and so $c'_4$ must be adjacent to one of $c_2$ or $c_3$, say $c_2$.
Then $\widehat{G}_{a_1,b'_2}$ has 9 edges and a 3-cycle $(c_2 b'_1 c'_4)$.
In the second case, if both $c'_4, c'_5$ are adjacent to another degree 3 vertex $c_2$,
then $\widehat{G}_{a_1,b'_1}$ has at most 9 edges and a 3-cycle $(c_2 c'_4 c'_5)$.
So we may assume that $c'_4 \adj \{c_2, c_3\}$ and $c'_5 \adj \{c_4, c_5\}$.
Again, we partially construct $\widehat{G}_{a_1,c'_4}$ as in Figure~\ref{fig:106025}(b).
Then $\widehat{G}_{a_1,c'_5}$ has 9 edges and a 3-cycle $(c_2 b'_1 c'_4)$.

Suppose instead that
$c'_1 \adj \{c_1, c_2\}$, $c'_2 \adj \{c_3, c_4\}$ and $c'_3 \adj \{c_5, c_6\}$.
If both $c_1, c_2$ are adjacent to $b'_1$ (similarly for $b'_2, c'_4 c'_5$),
then $\widehat{G}_{a_1,b'_2}$ has at most 9 edges and a 3-cycle $(c_1 c_2 b'_1)$.
Therefore we may assume that $b'_1 \adj \{c_1, c_3, c_5\}$.
Now if $b'_2$ (similarly for $c'_4, c'_5$) is adjacent to at least two among $c_1, c_3, c_5$,
say $c_1, c_3$, as in Figure~\ref{fig:106025}(c), 
then $\widehat{G}_{a_1,c_5}$ has 9 edges and a 3-cycle $(c_1 c_3 b'_2)$.
So we may assume that $c_1 \adj b'_2$, $c_3 \adj c'_4$ and $c_5 \adj c'_5$.
Now use the restoring method so that
we construct $G_{a_1,b'_1}$ as drawn in Figure~\ref{fig:106025}(d). 
After recovering $G$, we have an intrinsically knotted graph, 
from which we obtain the Heawood graph by deleting two edges connecting $a_1$ and $\{b'_1, b'_2\}$.

\section{$G$ only contains vertices with degree 3 or 4} \label{sec:g4}

In this section we assume that both $A$ and $B$ only contain vertices with degree 3 or 4.
Then $[A]$ and $[B]$ are either $[0,5,1]$ or $[0,2,5]$,
and we have three cases in the following subsections.
In Subsection~\ref{subsec:89} we find three IK graphs, each formed by adding two edges 
to the Heawood graph, see Figures~\ref{fig:025025}(f), (j), and (o).

\subsection{$[A]=[0,5,1]$ and $[B]=[0,5,1]$.}\ 

\begin{figure}[h!]
\includegraphics{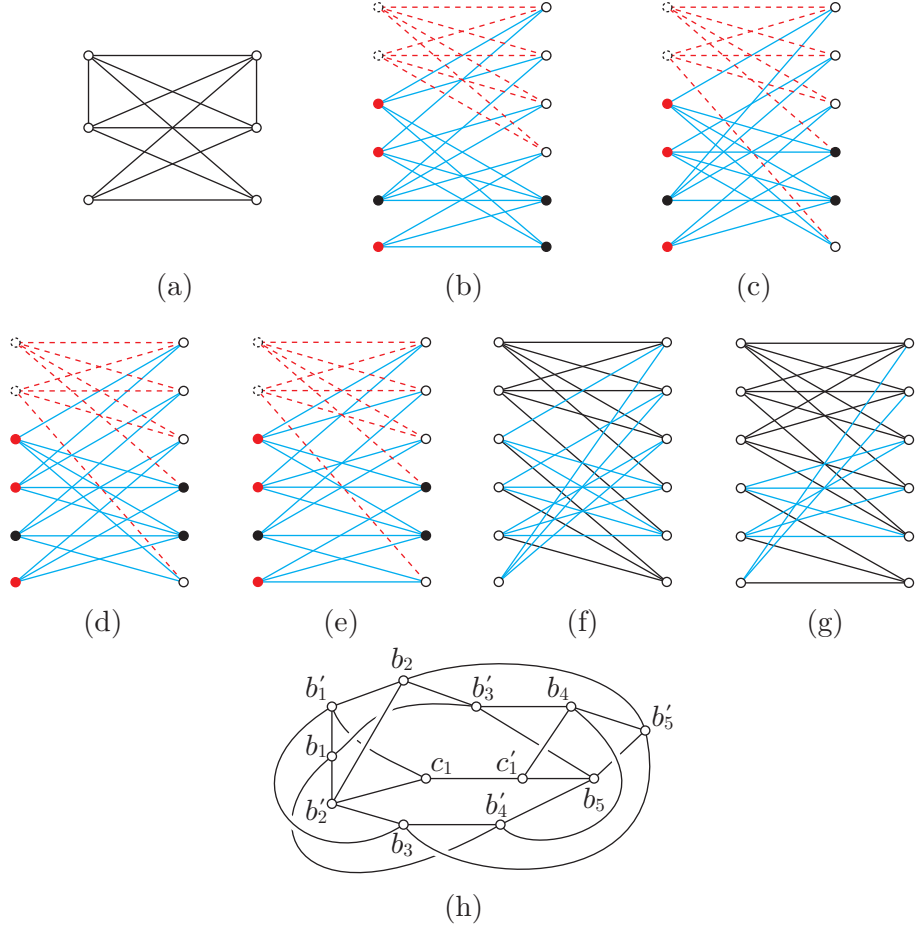}
\caption{Case of $[A]=[0,5,1]$ and $[B]=[0,5,1]$.}
\label{fig:051051}
\end{figure}

First we remark that if
$H$ is a graph, allowing multi-edges, 
that consists of four degree 4 vertices, two degree 3 vertices and eleven edges,
and such that there is a degree 4 vertex adjacent to the other three degree 4 vertices as well as a degree 3 vertex,
then the graph is non-planar only when it is the graph in Figure~\ref{fig:051051}(a).

Without loss of generality, there is a vertex $b_1$ so that 
$b_1 \adj \{b'_1, b'_2, b'_3, b'_4\}$.

As a first case, assume that some $b_2 \adj \{b'_1, b'_2, b'_3, b'_4\}$, and so $b'_5 \adj \{b_3, b_4, b_5, c_1\}$.
Using the restoring method, 
we construct $G_{b_1,b_2}$ so that $\widehat{G}_{b_1,b_2}$ is the $H$ mentioned in the remark above,
as shown in Figure~\ref{fig:051051}(b). 
After recovering $G$, we have a planar graph $\widehat{G}_{b_1,b_3}$.

As a second case, assume that $b_2 \adj \{b'_1, b'_2, b'_3, c'_1\}$, and so $b'_5 \adj \{b_3, b_4, b_5, c_1\}$ again.
Using the restoring method, 
we construct $G_{b_1,b_2}$ so that $\widehat{G}_{b_1,b_2}$ is the $H$ mentioned in the remark above,
as shown in Figure~\ref{fig:051051}(c), (d) and (e). 
After recovering $G$, we have planar graphs $\widehat{G}_{b_1,b_3}$ in all cases.

If $c_1 \nadj c'_1$, we assume $c'_1 \adj \{b_3,b_4,b_5\}$.
Then, to avoid the first case, we may say $b_2 \adj \{b'_1, b'_2, b'_3, b'_5\}$.
Furthermore, to avoid the second case, we conclude that 
$b'_4 \adj \{b_1, b_3, b_4, b_5\}$ and $b'_5 \adj \{b_2, b_3, b_4, b_5\}$.
By connecting the remaining edges, we obtain Figure~\ref{fig:051051}(f).
After recovering $G$, we have a planar graph $\widehat{G}_{b_1,b_3}$.

For otherwise, $c_1 \adj c'_1$, so may assume $c'_1 \adj \{b_4,b_5, c_1\}$, 
then to avoid the first case, we may say $b_2 \adj \{b'_1, b'_2, b'_3, b'_5\}$ and
$b_3 \adj \{b'_1, b'_2, b'_4, b'_5\}$.
Furthermore, to avoid the second case, we conclude that 
$\{b_4, b_5 \} \adj \{b'_3, b'_4, b'_5, c'_1\}$.
In a similar way, we obtain Figure~\ref{fig:051051}(g).
This graph has an embedding which does not have any non-trivially knotted cycles as shown in Figure~\ref{fig:051051}(h).

\subsection{$[A]=[0,5,1]$ and $[B]=[0,2,5]$.}\ 

First we remark that 
$NV_3(b_i,b_j) + NV_4(b_i,b_j) \neq 6$ for any pair $i,j$.
For otherwise, $\widehat{G}_{b_i,b_j}$ has 9 edges and at least two degree 4 vertices
because $V_3(b_i,b_j)$ has at most one vertex.

In this case, $b'_1$ and $b'_2$ are adjacent to at least two vertices in $A$ simultaneously.
First assume that there are exactly two common neightbors in $A$, one each of degree 3 and 4:
say, $b'_1 \adj \{b_1, b_2, b_3, c_1\}$ and $b'_2 \adj \{b_1, b_4, b_5, c_1\}$.
If $NV_3(b_2,b_3) = 5$, then $NV_3(b_2,b_3) + NV_4(b_2,b_3) = 6$.
Therefore $V_3(b_2,b_3)$ has at least two vertices;
say, $\{ b_2, b_3 \} \adj \{c'_1, c'_2 \}$ and similarly $\{ b_4, b_5 \} \adj \{c'_3, c'_4 \}$. 
We further assume that $b_1 \adj c'_1$.

\begin{figure}[h!]
\includegraphics{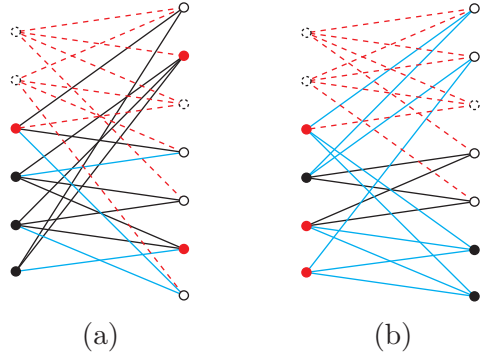}
\caption{Case of $[A]=[0,5,1]$ and $[B]=[0,2,5]$.}
\label{fig:051025}
\end{figure}

Now we distinguish three sub-cases.
First assume that $b_1 \adj c'_2$.
Then $b_4 \adj c'_5$, and so $NV_3(b_1,b_4) + NV_4(b_1,b_4) = 6$.
Second, assume that $b_1 \adj c'_3$ (similarly for $b_1 \adj c'_4$).
We may assume that $c'_5 \adj b_2$.
Using the restoring method,
we construct $G_{b_1,b_2}$ as shown in Figure~\ref{fig:051025}(a). 
After recovering $G$, we have a planar graph $\widehat{G}_{b_4,b_5}$.
Thirdly, we assume that $b_1 \adj c'_5$.
If $b_4 \nadj c'_5$, then $b_4 \adj c'_2$ and thus $NV_3(b_1,b_4) + NV_4(b_1,b_4) = 6$.
Therefore $b_4 \adj c'_5$, and similarly $b_5 \adj c'_5$.
By connecting the remaining edges, we obtain the same graph as Figure~\ref{fig:051025}(a)
after relabelling.

Without loss of generality, 
we assume that $\{ b'_1, b'_2 \} \adj \{b_1, b_2 \}$.
From the above remark, $NV_3(b_1,b_2) \leq 3$.
Now assume that $NV_3(b_1,b_2)=3$; say, $b_1 \adj \{ c'_1, c'_2 \}$ and $b_2 \adj \{ c'_1, c'_3 \}$.
If $NV_Y(b_1,b_2) = 1$, then again $\widehat{G}_{b_1,b_2}$ has 9 edges and a degree 4 vertex $b_3$.
So assume that $c'_1 \adj b_3$.
If $b_4 \nadj c'_2$, then we have $NV_3(b_1,b_4) + NV_4(b_1,b_4) = 6$.
Therefore $b_4 \adj c'_2$, and similarly $\{b_4, b_5\} \adj \{c'_2, c'_3\}$.
Using the restoring method,
we construct $G_{b_1,b_2}$ as shown in Figure~\ref{fig:051025}(b),
so it must be isomorphic to $K_{3,3}+e_2$.
Then $\widehat{G}_{b_3,b_5}$ has 9 edges and a degree 4 vertex $b_1$.

Finally we assume that $NV_3(b_1,b_2)=2$; say, $\{ b_1, b_2 \} \adj \{ c'_1, c'_2 \}$.
There is a degree 4 vertex $b_i$ among $b_3, b_4, b_5$ which is not adjacent to $c'_1, c'_2$.
Obviously, $NV_3(b_1,b_i) + NV_4(b_1,b_i) = 6$.

\subsection{$[A]=[0,2,5]$ and $[B]=[0,2,5]$, and the twin restoring method.}\ \label{subsec:89}

If there is a degree 4 vertex $b_i$ in $A$ which is not adjacent to degree 4 vertices in $B$,
or vice versa, then $NV_3(b_i,b'_1) \geq 7$, implying $|\widehat{E}_{b_i,b_1}| \leq 8$.
Therefore we assume that $b_1 \adj b'_1$ and $b_2 \adj b'_2$.
Now we distinguish three cases:
(1) $b_1 \nadj b'_2$ and $b_2 \nadj b'_1$,
(2) $b_1 \adj b'_2$ and $b_2 \nadj b'_1$, or
(3) $b_1 \adj b'_2$ and $b_2 \adj b'_1$.

\begin{figure}[h!]
\includegraphics{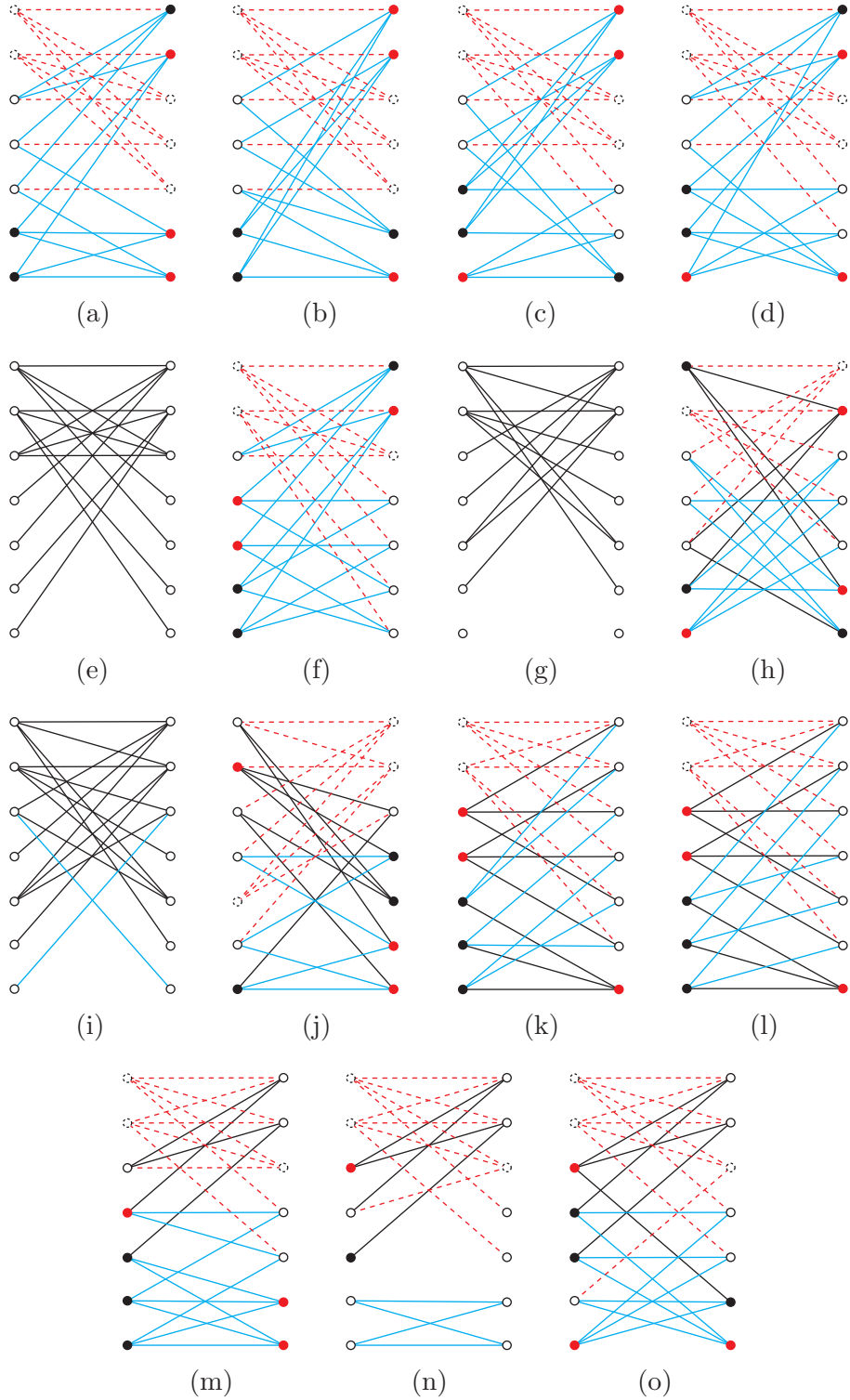}
\caption{Case of $[A]=[0,2,5]$ and $[B]=[0,2,5]$.}
\label{fig:025025}
\end{figure}

In case (1), 
we first assume that $NV_3(b_1,b_2) =3$, say $\{b_1, b_2\} \adj \{ c'_1, c'_2, c'_3\}$.
If some two vertices among $c'_1, c'_2, c'_3$ are adjacent to the same degree 3 vertex $c_i$ in $A$,
then $\widehat{G}_{b_1,b'_2}$ has 9 edges and a 2-cycle $(b_2 c_i)$.
Thus we further assume that $c_1 \adj c'_1$, $c_2 \adj c'_2$, and $c_3 \adj c'_3$.
Using the restoring method, we construct two graphs $G_{b_1,b_2}$ 
as shown in Figure~\ref{fig:025025}(a) and (b).
In both cases, after recovering $G$, we have planar graphs $\widehat{G}_{b'_1,c'_3}$.
We now assume that $NV_3(b_1,b_2) =4$, say $b_1 \adj \{ c'_1, c'_2, c'_3\}$ and
$b_2 \adj \{ c'_1, c'_2, c'_4\}$.
From the same argument above, we further assume that $c_1 \adj c'_1$ and $c_2 \adj c'_2$.
If $c_1 \adj c'_3$, then $\widehat{G}_{b_2,b'_1}$ has 9 edges and a 3-cycle $(b_1 c_1 c'_3)$.
Therefore $c_1 \nadj c'_3$, and similarly $c_1 \nadj c'_4$, $c_2 \nadj c'_3$, and $c_2 \nadj c'_4$.
Using the restoring method, we construct two graphs $G_{b_1,b_2}$ 
where there is no connection between degree 2 vertices in $A$ and $B$
as drawn in Figure~\ref{fig:025025}(c) and (d).
In both cases, after recovering $G$, we have planar graphs $\widehat{G}_{b'_1,c'_3}$.
Lastly we have $NV_3(b_1,b_2) = NV_3(b'_1,b'_2) =5$,
and we may have the connections shown in Figure~\ref{fig:025025}(e).
Here if $c_1 \nadj c'_1$, then there is no proper bipartition for $\widehat{G}_{b_1,b_2}$.
If $c_i \adj \{c'_2, c'_3\}$ (similarly for $\{c'_4, c'_5\}$), 
then $\widehat{G}_{b'_1,b'_2}$ has 9 edges and a 3-cycle $(b_1 c'_2 c'_3)$.
Thus all $c_i$, $i=2,3,4,5$, are adjacent to one of $\{c'_2, c'_3\}$ and another one of $\{c'_4, c'_5\}$.
Similarly all $c'_i$, $i=2,3,4,5$, are adjacent to one of $\{c_2, c_3\}$ and one of $\{c_4, c_5\}$.
Using the restoring method, we construct a graph $G_{b_1,b_2}$ satisfying the above conditions
as shown in Figure~\ref{fig:025025}(f).
After recovering $G$, we have an intrinsically knotted graph, 
from which we obtain the Heawood graph by deleting two edges:
one connecting $b_1$ and $b'_1$ and the other between $b_2$ and $b'_2$.

In the remaining cases (2) and (3), we remark that
both $\widehat{G}_{b_1,b_2}$ and $\widehat{G}_{b'_1,b'_2}$ have at most 9 edges.
It is sufficient to show that $NV_3(b_1,b_2) + NV_4(b_1,b_2) + NV_Y(b_1,b_2) \geq 6$
since $NE(b_1,b_2) = 8$.
If $NV_3(b_1,b_2) + NV_4(b_1,b_2) = 5$, then $NV_Y(b_1,b_2) = 1$, so we are done.
If $NV_3(b_1,b_2) + NV_4(b_1,b_2) = 4$, then  $V_3(b_1,b_2)$ has exactly two vertices, say $c'_1, c'_2$.
Let $c'_1 \adj c_1$.
If $c_1 \nadj c'_2$, then $NV_Y(b_1,b_2) = 2$.
For otherwise, when the remaining edge adjacent to $c_1$ is deleted in $\widehat{G}_{b_1,b_2}$, 
we can eventually delete one more edge.

In case (2), 
we first assume that $NV_3(b_1,b_2) =3$, there are two degree 3 vertices in $B$ 
such that $\{b_1, b_2\} \adj \{ c'_1, c'_2\}$.
Then $\widehat{G}_{b'_1,b'_2}$ has 9 edges by the remark and a 3-cycle $(b_2 c'_1 c'_2)$.

Now assume that $NV_3(b_1,b_2) =5$, say $b_1 \adj \{c'_1, c'_2\}$ and $b_2 \adj \{c'_3, c'_4, c'_5\}$.
By considering $\widehat{G}_{b_1,b_2}$, there are no degree 2 vertices in $A$.
Therefore we may assume that $b'_1, c_1$ and $c_2$ lie in the same partition,
implying $b'_1 \adj \{c_3, c_4, c_5\}$ as in Figure~\ref{fig:rest3}(b). 
The rest of process follows the twin restoring method, 
described in Subsection~\ref{subsec:restoring3} as an example.
Eventually we conclude that $\widehat{G}_{b'_1,b'_2}$ has 9 edges and a 5-cycle.

Finally we may assume that $NV_3(b_1,b_2) = NV_3(b'_1,b'_2) = 4$,
say $b_2 \adj \{c'_1, c'_2, c'_3\}$, $b_1 \adj \{c'_3, c'_4\}$,
$b'_1 \adj \{c_1, c_2, c_3\}$, and $b'_2 \adj \{c_3, c_4\}$ as in Figure~\ref{fig:025025}(g).
We remark that each of $c_1, c_2, c_4$ (similarly for $c'_1, c'_2, c'_4$) is adjacent to at most one of $c'_1, c'_2, c'_3$, 
for otherwise, $\widehat{G}_{b'_1,b'_2}$ has 9 edges and a 2 or 3-cycle.
If $c_3 \adj c'_3$ or $c_3 \adj c'_4$, then $\widehat{G}_{b_2,b'_1}$ has 9 edges and a 2 or 3-cycle.
Assume that $c_3 \adj c'_5$.
In this case, $b_1$ and $c_5$ are adjacent in $\widehat{G}_{b_2,b'_1}$,
and by the above remark, $c'_3 \adj c_5$ directly.
By applying the same remark again, $c_4 \adj c'_4$.
Without loss of generality, $c_4 \adj c'_1$ and $c'_4 \adj c_1$,
and thus $c_1 \adj c'_5$ and $c'_1 \adj c_5$.
By connecting the remaining edges, we obtain the graph as Figure~\ref{fig:025025}(h).
Then we have a planar graph $\widehat{G}_{b'_1,c'_3}$.
Therefore we may assume that $c_3 \adj c'_1$ and $c'_3 \adj c_1$ as Figure~\ref{fig:025025}(i).
Furthermore, we have $c_1 \adj c'_5$ because 
if $c_1 \adj c'_1, c'_2$ or $c'_4$, then $\widehat{G}_{b'_1,b'_2}$ has 9 edges and a 2 or 3-cycle.
Similarly we have $c'_1 \adj c_5$.
By the bipartition for $\widehat{G}_{b'_1,b'_2}$,
$c'_2 \adj \{ c_2, c_4\}$, and similarly $c_2 \adj c'_4$.
By connecting the remaining edges, we obtain the graph as Figure~\ref{fig:025025}(j).
After recovering $G$, we have an intrinsically knotted graph, 
from which we obtain the Heawood graph by deleting two edges:
one connecting $b_1$ and $b'_1$, the other connecting $b_2$ and $b'_2$.

In case (3),  
we first assume that $NV_3(b_1,b_2) =2$, there are two degree 3 vertices in $B$ 
such that $\{b_1, b_2\} \adj \{ c'_1, c'_2\}$.
Then $\widehat{G}_{b'_1,b'_2}$ has 9 edges by the remark and a 2-cycle $(c'_1 c'_2)$.

Now assume that $NV_3(b_1,b_2) =4$, say $b_1 \adj \{c'_1, c'_2\}$ and $b_2 \adj \{c'_3, c'_4\}$.
By considering $\widehat{G}_{b_1,b_2}$, there are no degree 2 vertices in $A$.
Therefore we may assume that $c_3,c_4,c_5$ lie in the same partition,
implying $c'_5 \adj \{c_3, c_4, c_5\}$. 
Now, we remark that
none of the five $c_i$'s can be adjacent to both $c'_1$ and $c'_2$ 
(and similarly for both $c'_3$ and $c'_4$),
for otherwise, $\widehat{G}_{b'_1,b'_2}$ has a 3-cycle $(c_i c'_1 c'_2)$.
Without loss of generality,
we say that $c_1 \adj \{b'_1, c'_1, c'_3\}$ and $c_2 \adj \{b'_2, c'_2, c'_4\}$.
We distinguish two cases:
$b'_1$ and $b'_2$ are adjacent to the same vertex $c_3$
or to different vertices, $b'_1 \adj c_3$ and $b'_2 \adj c_4$.
By connecting the remaining edges, we obtain the graphs in Figure~\ref{fig:025025}(k) and (l), respectively.
In both cases, after recovering $G$, we have planar graphs $\widehat{G}_{b_1,c_3}$.

Finally assume that $NV_3(b_1,b_2) = NV_3(b'_1,b'_2)  = 3$, 
say $b_1 \adj \{c'_1, c'_2\}$, $b_2 \adj \{c'_1, c'_3\}$, $b'_1 \adj \{c_1, c_2\}$, and $b'_2 \adj \{c_1, c_3\}$.
If $c'_1 \adj c_1$,
then using the restoring method, we construct a graph $G_{b_1,b_2}$ as drawn in Figure~\ref{fig:025025}(m).
After recovering $G$, we have planar graphs $\widehat{G}_{b'_1,b'_2}$,
which has 9 edges and a 2-cycle $(c'_2 c'_3)$.
If $c'_1 \adj c_2$, we partially construct $\widehat{G}_{b_1,b_2}$ as in Figure~\ref{fig:025025}(n).
As in the figure, $c'_4$ and $c'_5$ lie in the same partition, and so are $c_4$ and $c_5$.
Then $\widehat{G}_{b'_1,b'_2}$ has 9 edges and a 3-cycle;
in details, 
if $c_1$ and $c'_4$ lie in different partitions, then $c_1 \adj c'_4$ (or $c_1 \adj c'_5$),
in which case the 3-cycle is $(c_4 c_5 c'_5)$,
or otherwise, $c_3$ and $c'_4$ lie in different partitions,
then $c_3 \adj \{c'_4, c'_5\}$, in which case the 3-cycle is $(c_4 c'_4 c'_5)$.
Therefore we may assume that $c_1 \adj c'_4$ and $c'_1 \adj c_4$.
Using the restoring method, we construct a graph $G_{b_1,b_2}$ as drawn in Figure~\ref{fig:025025}(o).
After recovering $G$, we have an intrinsically knotted graph, 
in which we obtain $C_{14}$ by deleting two edges;
one connecting $b_1$ and $b'_1$, and another connecting $b_2$ and $b'_2$.


\begin{thebibliography}{99}
\bibitem{BM} J. Barsotti and T. Mattman,
    {\em Graphs on 21 edges that are not $2$--apex},
    Involve \textbf{9} (2016) 591--621.
\bibitem{BBFFHL} P. Blain, G. Bowlin, T. Fleming, J. Foisy, J. Hendricks and J. LaCombe,
    {\em Some results on intrinsically knotted graphs},
    J. Knot Theory Ramifications \textbf{16} (2007) 749--760.
\bibitem{FMMNN} E. Flapan, T. Mattman, B. Mellor, R. Naimi and R. Nikkuni,
    {\em Recent developments in spatial graph theory\/} 
    in Knots, links, spatial graphs, and algebraic invariants, 
    Contemp. Math. \textbf{689} (2017) 81--102.
\bibitem{GMN} N. Goldberg, T. Mattman and R. Naimi,
    {\em Many, many more intrinsically knotted graphs},
    Algebr. Geom. Topol. \textbf{14} (2014) 1801--1823.
\bibitem{HAMM} S. Huck, A. Appel, M-A. Manrique and T. Mattman,
    {\em A sufficient condition for intrinsic knotting of bipartite graphs}, 
    Kobe J. Math. \textbf{27} (2010) 47--57.
\bibitem{JKM} B. Johnson, M. Kidwell and T. Michael,
    {\em Intrinsically knotted graphs have at least $21$ edges},
    J. Knot Theory Ramifications \textbf{19} (2010) 1423--1429.
\bibitem{KLLMO} H. Kim, Lee, Lee, T. Mattman and S. Oh,
   {\em A new intrinsically knotted graph with 22 edges},
   Topology Appl. \textbf{228} (2017) 303--317.
\bibitem{KMO} H. Kim, T. Mattman and S. Oh,
    {\em Bipartite intrinsically knotted graphs with 22 edges},
    J. Graph Theory \textbf{85} (2017) 568--584.
\bibitem{KMO2} H. Kim, T. Mattman and S. Oh,
	{\em More intrinsically knotted graphs with 22 edges and the restoring method},
    J. Knot Theory Ramifications \textbf{27} (2018) 1850059.
\bibitem{KMO3} H. Kim, T. Mattman and S. Oh,
   {\em Triangle-free intrinsically knotted graphs with 22 edges and the twin-restoring method},
   (in preparation)
\bibitem{LKLO} M. Lee, H. Kim, H. J. Lee and S. Oh,
    {\em Exactly fourteen intrinsically knotted graphs have 21 edges},
    Algebr. Geom. Topol. \textbf{15} (2015) 3305--3322.
\bibitem{M} T. Mattman,
    {\em Graphs of 20 edges are 2-apex, hence unknotted},
    Algebr. Geom. Topol. \textbf{11} (2011) 691--718.
\bibitem{OT} M. Ozawa and Y. Tsutsumi,
    {\em Primitive spatial graphs and graph minors},
    Rev. Mat. Complut. \textbf{20} (2007) 391--406.
\bibitem{RS} N. Robertson and P. Seymour,
    {\em Graph minors XX, Wagner's conjecture},
    J. Combin. Theory Ser. B \textbf{92} (2004) 325--357.
\end{thebibliography}
\end{document}